\documentclass{amsart}
  \usepackage{amscd,amssymb,epsfig}
   \usepackage{epic,eepic}

                     \numberwithin{equation}{subsection}

                     \newtheorem{propo}{Proposition}[subsection]
                     \newtheorem{corol}[propo]{Corollary}
                     \newtheorem{theor}[propo]{Theorem}
                     \newtheorem{lemma}[propo]{Lemma}
                     \theoremstyle{definition}

                     \theoremstyle{remark}

		     \newcommand{\CC}{\mathbb{C}}
		     \newcommand{\QQ}{\mathbb{Q}}
                     \newcommand{\ZZ}{\mathbb{Z}}
                     \newcommand{\RR}{\mathbb{R}}

		     \newcommand{\A}{\mathcal{A}}
 \newcommand{\M}{\mathcal{M}}
		     \newcommand{\str}{\mathcal{S}}

                     \newcommand{\Hom}{\operatorname{Hom}}
		     \newcommand{\Ker}{\operatorname{Ker}}
		     \newcommand{\Exp}{\operatorname{Exp}}
\newcommand{\Lbl}{\operatorname{Lbl}}
\newcommand{\lbl}{\operatorname{lbl}}

                     \newcommand{\id}{\operatorname{id}}

 \newcommand{\rrr}{r}

		     \newcommand{\sign} {\operatorname {sign}}
		    \newcommand{\modu}{\operatorname{mod}}
\newcommand{\Int}{\operatorname{Int}}
		     \newcommand{\arr}{\operatorname{arr}}
		  \newcommand{\rank}{\operatorname{rank}}
 \newcommand{\edg}{\operatorname{edg}}
		  \newcommand{\Perm}{\operatorname{Perm}}

\begin{document}
      \title{Virtual strings}
                     \author[Vladimir Turaev]{Vladimir Turaev}
                     \address{%
              IRMA, Universit\'e Louis  Pasteur - C.N.R.S., \newline
\indent  7 rue Ren\'e Descartes \newline
                     \indent F-67084 Strasbourg \newline
                     \indent France }
                     \begin{abstract} 	A virtual string    is   a scheme of  self-intersections of a closed
curve on  a  surface.   We  introduce virtual strings and study their  geometric properties and homotopy invariants. 
 We also discuss connections  between virtual strings, Gauss words, and virtual
knots. 
                     \end{abstract}
                     \maketitle

\centerline {\bf Contents}
 \vskip1truecm

\noindent  {\bf  1. Introduction}

\noindent {\bf 2.  Generalities on virtual strings}

\noindent  {\bf 3.  Polynomial $u$}

\noindent {\bf 4.  Geometric realization  of virtual strings}

\noindent {\bf 5.  Combinatorics of closed curves on the 2-sphere}

\noindent {\bf 6.  Based  skew-symmetric  matrices}

\noindent  {\bf 7.  Based   matrices  of  strings}

\noindent {\bf 8.  Lie cobracket for strings}

\noindent {\bf 9.  Virtual strings versus  virtual  knots}

\noindent {\bf 10.  Proof of Theorem \ref{th:eee}}

\noindent {\bf 11.   Algebras and groups associated with strings}

\noindent {\bf 12.   Open questions}

                  \section{Introduction} 
		  
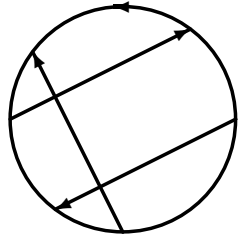
\begin{figure}[ht]
                     \setlength{\unitlength}{0.5cm}
                     \begin{picture}(18,6)(-10,-2.5)
                             \allinethickness{.5mm}
                     %Circle
                             \put(0,0){\circle{6}}
                             \put(0,3){\vector(-1,0){.25}}
                             \put(-3,0){\vector(2,1){4.8}}
                             \put(0,-3){\vector(-1,2){2.4}}
                             \put(3,0){\vector(-2,-1){4.8}}
                     \end{picture}
                     \caption{A virtual string of rank 3}\label{fg:gdg}
                     \end{figure}

	A virtual string    is   a scheme of  self-intersections of a generic oriented closed
curve on  an 
oriented surface. More precisely, a    virtual string of rank $m\geq 0$   
is  an 
oriented circle  with $2m$ distinguished points  partitioned into $m$ ordered pairs. These 
$m$ ordered pairs   of points  are 
called   arrows of the virtual string. An example of a virtual string of rank 
$3$  is 
shown on Figure~\ref{fg:gdg} where the  arrows are represented by 
geometric 
vectors.
 
A  (generic  oriented)  closed curve  on an oriented 
surface   gives rise to an \lq\lq underlying" virtual string whose arrows correspond to the self-crossings of the curve.  The
usual homotopy of curves   suggests a  notion of homotopy for   strings.  The main objective of the theory of virtual
strings is a study (and  eventually classification) of their homotopy classes. To this end, we  
introduce    certain   homotopy invariants of virtual 
strings, specifically a one-variable polynomial  $u$ and a so-called
based matrix.   

In analogy with the theory of knots we  define a class of slice virtual strings.  A string is slice if it
underlies a closed curve on a closed surface  $\Sigma$ which is contractible in a handlebody bounded by $\Sigma$. We
formulate  obstructions to the sliceness of  a string in terms of the polynomial $u$ and the based matrix. 

We introduce a natural Lie cobracket in  the free abelian group generated by the homotopy classes of strings. Dually, the
abelian group  of   $\ZZ$-valued homotopy invariants of strings  becomes
 a Lie algebra.    This  Lie algebra gives rise to  further  algebraic objects including a Hopf algebra structure in the
(commutative) polynomial algebra generated by the homotopy classes of strings.

		     	    The notion of a virtual  string  should be compared with the  one        of  a Gauss word. A  Gauss word (or a
Gauss code)  is a 
sequence of letters in a finite alphabet   in which all letters of the 
alphabet 
occur  exactly twice. The sequence is considered up to circular 
permutations.  Gauss    observed that a    closed curve on $\RR^2$    gives rise   to a Gauss word. 
Namely, label the  self-crossings of 
the curve by different letters and write down  these letters in the 
order of 
their appearance when one traverses   the whole curve   beginning at  a 
generic 
point.   The 
resulting word   is well defined  up to circular permutations (and the choice 
of 
letters). A similar procedure applies to virtual strings  where instead of crossings one should label the arrows. 
Thus, each virtual string gives rise to a Gauss word. In fact,  the notion of a Gauss word is equivalent to the 
notion of 
a virtual string with the order in  the distinguished pairs  of  points 
forgotten. 

Gauss  introduced his sequences of letters  in an attempt to give a   combinatorial
formulation of closed curves on $\RR^2  $.  We discuss another relevant ingredient: a bipartition of the alphabet.
As an application of the theory of virtual strings, we obtain  a complete  combinatorial description of closed curves on
the 2-sphere   in terms of   Gauss words and   bipartitions. This description extends the   theorem of
P. Rosenstiehl
\cite{ro}   
  characterizing  the  Gauss words realizable by closed curves on $\RR^2$.

 Virtual strings are   closely related to   virtual knots   introduced by L.  Kauffman \cite{Ka}. 
		 They can be defined in terms of  so-called arrow diagrams which are just   virtual strings whose arrows  are 
provided with signs $+$ or $-$.  Virtual knots are arrow diagrams considered up to several moves, induced by the
Reidemeister moves on ordinary knot diagrams. The  term  \lq\lq virtual knots" 
suggested to us the term virtual strings.  

Forgetting the signs of the arrows  transforms an arrow diagram into  a virtual string.  This gives a map from the set of
virtual knots into the set of  homotopy classes of virtual strings. We give a more elaborate construction which associates
with each virtual knot  a polynomial  expression in virtual strings with coefficients in the ring
$\QQ[z]$. This leads to  an isomorphism between  a skein algebra of virtual knots, defined below,  and  a polynomial algebra
generated by the homotopy classes of strings. 
 
 The organization of the paper should be clear from the Contents above.
	
		       \section{Generalities on virtual strings}
                    
		         \subsection{Definitions}\label{sn:g11} We give 
here a formal definition of a virtual string.   For an  integer $m\geq 
0$,   a 
{\it  virtual string $\alpha$ of rank $m$} (or briefly a {\it 
string}) 
is  an oriented circle, $S$, called the {\it core circle} of $\alpha$,  and a 
distinguished   set of $2m$ distinct points of $S$ partitioned into $m$  
ordered 
pairs.    We   call these $m$ ordered     pairs   of points  the {\it 
arrows} of 
$\alpha$. The set of arrows of $\alpha$ is denoted $\arr(\alpha)$. The 
endpoints 
$a,b\in S$ of  an arrow  $(a, b)\in \arr(\alpha)$  are called   its 
{\it  tail} 
and {\it   head}, respectively. The $2m$ distinguished 	points of $S$ 
are called 
the {\it endpoints} of $\alpha$.	
			 
 The  string formed by an oriented circle and an empty set of 
arrows is 
called a {\it trivial virtual string}.  
  An example of a virtual string of rank $3$  is shown on 
Figure~\ref{fg:gdg}.

			       By a {\it homeomorphism} of two virtual strings, 
we mean  an orientation-preserving  homeomorphism of the core circles 
transforming the  set of  arrows of the  first   string   onto the  set 
of arrows 
of the second string. Two virtual strings are {\it homeomorphic} if 
they are 
related  by   a homeomorphism. Clearly, homeomorphic strings have the same rank. 

By abuse of language, the homeomorphism classes of       virtual 
strings will be 
also called virtual strings.

		      \subsection{From curves   to    strings}\label{sn:g12} 
		     By a   surface,  we mean a smooth oriented 
2-dimensional manifold.   By a    {\it   closed curve} on a 
surface 
$\Sigma$, we mean a generic  smooth immersion  $\omega$ of an oriented circle $S$ into 
$\Sigma$. Recall that   a smooth map $ S\to \Sigma$ is  an {\it immersion} if   
its 
differential is  non-zero at all points of $S$.  An immersion 
$\omega:S\to \Sigma$ 
is {\it generic} if  $ \# (\omega^{-1}(x))\leq 2$  for all  $ x\in \Sigma$,  the 
set $\{x\in \Sigma\,\vert\, \# (\omega^{-1} (x))=2\}$ is finite, and   all 
its  points 
  are transverse intersections of   two branches.  Here and below the symbol 
$\# (A)$ 
denotes the cardinality of a set $A$. The points $x\in \Sigma$ such 
that $ \# 
(\omega^{-1} (x))=2 $ are called {\it double points} or {\it crossings} of 
$\omega$.  
		     
		     A closed curve $\omega:S\to \Sigma$   gives rise to an {\it underlying  virtual 
string} $\alpha_\omega$. The core circle of $\alpha_\omega$  is $S$ and the   arrows  of $\alpha_\omega$ 
are 
all ordered pairs $a,b\in S$ such that $\omega(a)=\omega(b)$ and the pair (a 
positive 
tangent vector of $\omega$ at $a$, a positive tangent vector of $\omega$ at 
$b$)  is a 
positive basis in   the tangent space of   $\omega(a)$.  For instance,  the underlying string of a simple closed curve on
$\Sigma$ is a trivial virtual string.

We say that a  virtual string  
 is {\it realized} by a closed curve $\omega:S\to \Sigma$ if  it is homeomorphic to $\alpha_\omega$. As we 
shall see 
below, every virtual string can be realized by a closed curve  on a 
surface.
		    
		     \subsection{Homotopy of   strings}\label{sn:g13}    The 
usual 
homotopy of closed curves on a surface suggests  to introduce a relation of homotopy for   virtual strings. Observe first   
that   two 
homotopic   curves on a surface can be related by a finite sequence of 
the 
following \lq\lq elementary" moves     (and the inverse 
moves):

		     (a) a local move adding a small curl to the curve;
		     
		     (b) a local move pushing a branch of the curve across 
another branch and creating  two new double points;
		     
		     (c) a local   move pushing a branch of the curve across a 
double point;
		     
		     (d) ambient isotopy in the surface.

	The move (a) has two forms (a${)}^{+}$   and (a${)}^{-}$   depending on	whether the curl lies on 
the left or   the right of  the curve where the left and the right are determined by the direction of the curve and the
orientation of the surface.  Considered up to ambient isotopy, the move (b) has three forms  depending on the direction 
of the two branches.   Similarly,  considered up to ambient isotopy, the move (c) has two
forms  (c${)}^{+}$   and (c${)}^{-}$  
depending on the direction of the   branches. Using the standard braid generators $\sigma_1, \sigma_2$ on 3 strands we
can encode  this move  as
$\sigma_1\sigma_2 
\sigma_1\mapsto \sigma_2\sigma_1 \sigma_2$ where the over/undercrossing 
information is forgotten.  The  moves  (c${)}^{+}$   and (c${)}^{-}$    are obtained by directing (before
and after the move)  the first and third strands up and the second strand up or down, respectively.   It is easy to see that  
(c${)}^{+}$,  (c${)}^{-}$  can be obtained from each other  using    ambient isotopy,    moves  (b), and  inverses to
 (b).  Similarly, the moves (a${)}^{+}$,  (a${)}^{-}$  can be obtained from each other  using    ambient isotopy,   
moves  (b),  (c${)}^{-}$, and  inverses to   (b).  Thus the moves (a${)}^{-}$, (b), (c${)}^{-}$ generate all the other moves.

It is clear that  ambient isotopy of a closed curve does not change 
the 
underlying virtual string. We now describe the analogues  for virtual 
strings of 
the moves (a${)}^{-}$, (b), (c${)}^{-}$. In this description and in the sequel, by an 
{\it arc} 
on  an oriented  circle $S$ we mean an {\it embedded arc} on $S$. The orientation 
of $S$ 
induces an orientation of all arcs on $S$. For two distinct points 
$a,b\in S$, we 
write   $ab$ for the unique oriented arc  in $S$ which      begins in 
$a$ and 
terminates in $b$.  Clearly, $S= ab \cup ba$ and $ab\cap ba= \{a,b\}$.

Let $\alpha$ be a virtual string with core circle $S$. Pick two 
distinct points 
$a,b \in S$  such that the arc   $ab\subset S$ is disjoint from the set of   
endpoints of 
$\alpha$. The move (a$)_s$, where $s$ stands for \lq\lq string",   
adds to 
$\alpha$  the pair $(a,b)$.    This amounts 
to 
attaching a small arrow to $S$  such that  the arc  in $S$ leading from its tail to its head  is disjoint from the   endpoints of
$\alpha$. The move (b$)_s$  acts on 
$\alpha$ as follows. Pick two   arcs   on $S$ disjoint from each 
other 
and   from the   endpoints of $\alpha$. Let $a,a'$ be the 
endpoints of 
the first arc (in an arbitrary order) and $b,b'$ be the endpoints of 
the second 
arc. The move   (b$)_s$ adds 
to $\alpha$ two arrows $(a,b)$ and $(b',a')$.  (This move has four forms depending on the
  two possible choices for  $a$  and two possible choices for  $b$.   However, two of these
 forms  of (b$)_s$  are  equivalent.) The move (c$)_s$     applies  to 
$\alpha$   when $\alpha$ has three arrows 
$(a^+,b), (b^+,c), (c^+,a)$ where  $a, a^+, b,b^+,c, c^+\in S$ such that
the   arcs   $aa^+$,   $bb^+$,     $cc^+$ are disjoint from each 
other and 
from the other   endpoints of $\alpha$. The   move (c$)_s$  replaces 
the  arrows  $(a^+,b), (b^+,c), (c^+,a)$ with the arrows  $ (a,b^+), (b,c^+), (c,a^+)$.

	We say that two virtual strings are {\it homotopic} if they can be 
related by a finite sequence of homeomorphisms, the  {\it homotopy moves} (a$)_s$, 
(b$)_s$, 
(c$)_s$, and the inverse moves. A virtual string homotopic to a trivial 
virtual 
string is said to be {\it homotopically trivial}. 

It is clear from what 
was said 
above that the underlying virtual strings of homotopic closed curves on 
a 
surface are themselves homotopic.  
	
\subsection{Transformations of  
strings}\label{sn:g1335}   For a   string $\alpha$, we define the {\it  opposite string}  $  
\alpha^-$ to be $\alpha$ with opposite orientation on 
the core 
circle. The {\it inverse  string} $ \overline \alpha$ is obtained from 
$\alpha$ 
by reversing all its arrows. On the level of closed curves on surfaces, 
these 
two transformations correspond to traversing the same curve in the 
opposite 
direction and to inverting the orientation of the ambient surface, respectively.
If two strings are homotopic, then their opposite (resp. inverse) 
strings are 
homotopic.  

One can raise a number of  questions concerning the transformations $\alpha\mapsto \alpha^-$,  $\alpha\mapsto \overline
\alpha$, $ \alpha\mapsto \overline
\alpha^-$. For instance, one can ask 
  whether  there is a  string $\alpha$ that is not homotopic to $\alpha^-$ (resp.  to $ \overline
\alpha$,  $\overline
\alpha^-$).    Below we will answer this
question   in the positive. 

A virtual string $\alpha$ with core circle $S$ is a {\it product} 
of 
virtual strings $\alpha_1$ and $\alpha_2$ if there are   disjoint 
arcs 
$a_1b_1, a_2b_2\subset S$ such that each arrow of $\alpha$ has both 
endpoints on either $a_1b_1$ or on $a_2b_2$ and the   string formed by 
$S$ and 
the arrows of $\alpha$ with endpoints on $ a_i b_i$ is homeomorphic to 
$\alpha_i$ for $i=1,2$.   
One can ask whether the product is  a  well-defined operation on strings (at least up to homotopy) and whether it is
commutative. Below we will answer these
questions  in the negative.

\subsection{Geometric invariants   of  
strings}\label{sn:g133} 
 We define four geometric characteristics of  strings: the genus, the 
homotopy 
genus,  the slice  genus, and the homotopy rank. For a string $\alpha$, its {\it genus} 
$g(\alpha)$  is  the minimal integer $g\geq 0$   such that $\alpha$ can 
be 
realized by a closed curve on a surface of genus $g$. The {\it homotopy 
genus} 
$hg(\alpha)$  is     the minimal integer $g\geq 0$   such that $\alpha$ 
is 
homotopic to a string of genus $g$.      The {\it homotopy rank} 
$hr(\alpha)$  
is  the minimal integer $m\geq 0$   such that $\alpha$ is homotopic to 
a string 
of rank $m$. For example, if $\alpha$ is a trivial string, then 
$g(\alpha)=hg(\alpha)=hr(\alpha)=0$.  
 It is clear that  the homotopy genus and the homotopy rank are 
homotopy 
invariants of strings. Below we  compute the genus  explicitly  and show that it is not a 
homotopy 
invariant. 
  
For $k\geq 0$ denote by $\Sigma_k$  a compact oriented surface  of genus $k$   bounded by a circle.  The {\it
slice  genus}
$sg(\alpha)$ of a string
$\alpha$ is the minimal  integer
$k\geq 0$ such that  there is a handlebody $H$ (of a certain genus) and a map $\Omega:\Sigma_k\to H$ such that
$\Omega(\partial \Sigma_k)\subset \partial H$ and the  map $\Omega\vert_{\partial \Sigma_k }:\partial \Sigma_k \to
\partial H$ is a (generic) closed curve on $\partial H$ representing $\alpha$.  The existence of such  $k$ follows
from the fact that   any loop on a closed surface becomes homologically trivial in a certain  handlebody bounded by this
surface.

 If $sg(\alpha)=0$, then we say that   $\alpha$ is {\it slice}. 
Thus,  a string is slice  if it can be realized on a closed surface by a  
curve   contractible in a handlebody bounded by this surface.  For example, a trivial
string is slice.

The genus, the homotopy genus,  the slice genus, and the  homotopy rank 
of a string are  preserved under the transformations $\alpha \mapsto \alpha^-,  \alpha \mapsto  \overline \alpha$.

 \subsection{Encoding of strings}\label{sn:g1345} There are two 
simple 
methods allowing to encode virtual strings in a compact way. Although 
we do 
not use these methods in this paper, we briefly describe them   for 
completeness.
 
(1)  Consider a 
finite set 
$E$ consisting of $m$ elements and its disjoint copy $E^+=\{x^+\,\vert 
\, x\in 
E\}$. 
Let $y_1,y_2,\ldots, y_{2m}$  be  a sequence of elements of the set 
$E\cup E^+$ 
in which every element  appears  exactly once. (Such a sequence 
determines a 
total order in  $E\cup E^+$ and vice versa.)  The sequence  $ y_1,y_2,\ldots, y_{2m}$ defines  a 
string  
of rank $m$ whose  underlying circle   is    $S=\RR\cup \{\infty\}$ 
with    
right-handed orientation on $\RR$ and whose  $m$   arrows are the pairs 
$(a,b)$ 
such that $a,b\in \{1,2,\ldots, 2m\}\subset S, y_a\in E$, and 
$y_b=y_a^+\in E^+$.  Any string can be encoded in this way. For instance, the 
string drawn 
in Figure \ref{fg:gdg} is encoded by the sequence
$x^+,y,z^+,x,z,y^+$ where $E=\{x,y,z\}$.

(2) By Section \ref{sn:g12},  virtual strings can 
be 
encoded by  closed curves on surfaces. This   has an extension similar 
to 
Kauffman's  graphical  encoding of virtual knots in \cite{Ka}. Namely, 
consider 
a (generic) closed 
curve   on a surface and suppose that some of its crossings are marked 
as \lq\lq 
virtual". Take the   string of this curve as in Section 
\ref{sn:g12} 
and forget all its arrows corresponding to virtual crossings. 
 It is easy to see that every virtual string can be obtained in this 
way from a  
 closed curve  in $\RR^2$ with virtual crossings. This  yields a  
graphical 
encoding of strings by plane curves with virtual crossings. 
The relation of  homotopy for strings has a simple description in  this 
language: it is generated by  the moves shown in \cite{Ka}, 
Figure 2
(where the over/undercrossing  information should be forgotten).

 \subsection{Remarks}\label{sn:g14}   1.    We can point out certain classes of closed  curves on surfaces whose 
underlying  virtual strings are homotopically trivial.   Since  all closed curves on  
$ S^2$ are contractible,  their underlying   strings are  homotopically trivial.  Therefore the same is true for 
closed  curves on any subsurface of $S^2$, i.e.,   on any surface of genus 0.  In particular, all  closed  curves  on an
annulus have   homotopically trivial underlying  strings.  Since each  closed curve  on a torus can be deformed into an
annulus,  its  underlying string   is homotopically trivial.  The
same holds for    closed curves on a torus with holes. 
            
 2.  The move (a$)_s$     has a   version  (a${)}^{+}_s$ which 
is defined   as (a$)_s$ above but adds the arrow $(b,a)$ rather than $(a,b)$.  
 This  move 
underlies 
the move (a${)}^{+}$ on closed curves.  The
move  (a${)}^{+}_s$ preserves the homotopy class of  a string.  Indeed, it     can be 
expressed  as a composition of   (b$)_s$,  (c$)_s$, and an inverse  to   (a$)_s$.

3. The move (c$)_s$     has a  version  (c${)}^{+}_s$ which 
applies  to a string    
 when it has three arrows 
$(a,b), (a^+,c)$, $(b^+,c^+)$ 
such that
the   arcs   $aa^+$,   $bb^+$,    $cc^+$ are disjoint from each 
other and 
from the other   endpoints of the string. The   move   (c${)}^{+}_s$  replaces these 
three arrows 
with the arrows  $ (a^+,b^+), (a,c^+), (b,c)$.
 This  move 
underlies 
the move (c${)}^{+}$ on closed curves.  The   move   (c${)}^{+}_s$  can be 
expressed
as a composition of   (c${)}^{-}_s$    and  (b$)_s$.
 
4. For any  virtual string $\alpha$, its appropriate product with $\overline \alpha^-$ is
slice.  Indeed, let  
$S$ be the core circle of $\alpha$  and let $ab\subset S$  be an arc containing all the   
endpoints of 
$\alpha$. Let $(\alpha', S', a'b'\subset S')$ be a disjoint copy of the triple $(\alpha, S, ab)$.  Consider the circle  
$S''= (ab\cup a'b')/a=a', b=b'$  and provide  it with the orientation extending the one on $ab$. The arrows
of
$\alpha$ and $\overline \alpha'$ are attached to $ab\cup a'b'$ and form in this way a virtual string, $\alpha''$, with core
circle 
 	$S''$.  It is clear that $\alpha''$ is a product of $\alpha$ with $\overline \alpha^-$. We claim that $\alpha''$ is slice.  To
see this, represent $\alpha$ by a closed curve $\omega:S\to \Sigma$ on a  surface $\Sigma$.  The map $\omega$
transforms   $S- ab$ onto an embedded arc in $\Sigma$ disjoint from the rest of the curve. Let $D\subset
\Sigma$ be a 2-disc such that  $D\cap \omega(S)=\omega(S-ab)$ and   $\partial D\cap
\omega(S)=\{\omega( a), \omega(b)\}$. The 3-manifold $H=(\Sigma-\Int D) \times [0,1]$ is a handlebody. The four paths
$\omega(ab)\times 0, \omega(ab)\times 1, \omega(a) \times [0,1], \omega(b) \times [0,1]$ form  a closed curve on
$\partial H$   realizing $\alpha''$ and contactible in $H$.

5. Replacing the circle  in the definition of a string  by an oriented one-dimensional manifold $X$ we obtain a  
{\it virtual string with core manifold} $X$.   The definitions and results of  this
paper can be extended to such strings  with appropriate changes.  Of special interest are strings with core manifold
 $X=[0,1]$;  we call them  {\it open strings}. In this context it is natural to call     virtual strings with core manifold a circle   
 {\it closed} strings.    Gluing the endpoints of 
$X=[0,1]$, we can transform an open string into a closed   string, its {\it closure}.  Open strings  underlie  (generic)
paths on  a surface with   endpoints on the boundary. Open strings  can be  multiplied in the obvious way.

 \subsection{Exercise}\label{sn:g1zz4}  Verify  that all virtual strings of rank 
$\leq 2$ 
are homotopically trivial.

            \section{Polynomial $u$}\label{abric} 
                    
		         \subsection{Invariants $\{u_k\}_k$}\label{sn:g21} Let 
$\alpha$ be a virtual string with core circle $S$. Each arrow 
$e=(a,b)\in 
\arr(\alpha)$ splits $S$ into two arcs   $ {ab}$ and $  {ba}$. We say 
that an 
arrow  $f=(c,d)$  of $\alpha$ (distinct from $e$) {\it links} $e$ if one of its endpoints lies on 
$ab$ and the 
other one lies on $ba$. More precisely,     $f=(c,d)$ 
  links $e$ {\it positively} (resp. {\it negatively}) if  $c\in 
ab, d\in 
ba$ (respectively, if $ c\in ba, d\in ab$).  If $f$ does not link $e$, then  $e$ and $f$ are {\it unlinked}.  Let
$n(e)\in
\ZZ$ be the algebraic  number of arrows of
$\alpha$  linking
$e$, i.e., 
  the 
number of arrows of $\alpha$  linking $e$ positively minus   the 
number of arrows of $\alpha$  linking $e$
negatively.
			 
			 It is easy to trace the behaviour of   $n(e)$ under the 
homotopy moves (a$)_s$, (b$)_s$,  (c$)_s$ on $\alpha$.
			 The move (a$)_s$ adds an arrow $e_0$ with  $ n( e_0)=0$ 
and keeps $n(e)$ for all other arrows. The move (b$)_s$ adds two arrows 
$e_1,e_2$ with $n(e_1)=-n( e_2)$ and keeps $n(e)$ for all other arrows. 
Consider 
  the move (c$)_s$ and use the notation of Section \ref{sn:g13}. It is 
obvious 
that for all arrows $e$ preserved under the move, the number $n(e)$ is 
also 
preserved.  
Each arrow $e=(a^+,b), (b^+,c), (c^+,a)$ occuring before the move gives 
rise to 
an arrow			 $ e'= (a,b^+), (b,c^+), (c,a^+)$, respectively, 
occuring after the move. We claim that $n(e)	=n(e')$.	Consider for 
concreteness $e=(a^+,b)$. Note that  the points $c,c^+$ lie either on 
$ab$ or on 
$ba$. Suppose that  $c,c^+\in ab$. Then  the arrows	$(b^+,c)$ and $ 
(c^+,a)$ 
contribute $1$ and $-1$ to $n(e)$, respectively, while the  corresponding 
arrows 
$(b,c^+)$ and $ ( c, a^+)$		contribute  $0$ to $n(e') $.
All other arrows contribute the same to  	$n(e)$ and  $n(e') $. 
Hence 
$n(e)=n( e')$. 
If $c,c^+\in ba$, then  the arrows	$(b^+,c)$ and $ (c^+,a)$ 
contribute $0$ to $n(e)$ while the corresponding arrows 
$(b,c^+)$ and $ ( c, a^+)$ 	contribute $-1$ and $1$ to $n( e')$, respectively.
All other arrows contribute the same to  	$n(e)$ and  $n( e')$. 
Hence 
$n(e)=n( e')$.  

For an integer $k\geq 1$,  set  $$u_k(\alpha)=\#\{e\in 
\arr(\alpha)\,\vert \,n(e) =k\}
- \#\{e\in \arr(\alpha)\,\vert \,n(e) =- k\} \in \ZZ.$$
  It is clear from what was said above that $u_k(\alpha)$ is preserved 
under the 
moves  (a$)_s$, (b$)_s$,  (c$)_s$. In other words, $u_k(\alpha)$ is a homotopy 
invariant 
of $\alpha$.  Clearly, $u_k(\alpha)=0$ for all $k$ greater than or equal to the rank of $\alpha$.  If $\alpha$ is homotopically
trivial,  then 
$u_k(\alpha)=0$ for all $k\geq 1$.

\subsection{Polynomial  $u(\alpha)$}\label{sn:g22} We can combine the     
invariants 
$u_k$ of a virtual string $\alpha$ into a polynomial 
$$u(\alpha)=\sum_{k\geq 1} u_k(\alpha)\, t^k$$
where $t$ is a  variable. The free term of this polynomial is always 
$0$ and its 
degree is bounded from above by $m-1$ where $m$ is the rank of 
$\alpha$. This 
polynomial is a homotopy invariant of $\alpha$. If $\alpha$ is 
homotopically 
trivial, then $u(\alpha)=0$. (The converse is not true, as we shall see 
below.) 
The polynomial $u(\alpha)$ yields an  estimate  for the 
homotopy rank 
$hr(\alpha)$ of $\alpha$ defined in Section \ref{sn:g133}: 
\begin {equation}\label{homr}
hr(\alpha)\geq \deg u(\alpha)+1.\end{equation}

We can rewrite  $u(\alpha)$ as follows:
\begin {equation}\label{homr11} u(\alpha)=\sum_{e\in \arr(\alpha), n(e)\neq 0} \sign (n(e))\, 
t^{\vert 
n(e)\vert}\end{equation}
where $\sign (n)=1$ for positive $n\in \ZZ$ and $\sign (n)=-1$ for 
negative 
$n\in \ZZ$.
 Therefore
$$\sum_{k\geq 1} k\, u_k(\alpha) \,t^{k-1}=u'(\alpha)=\sum_{e\in \arr(\alpha), n(e)\neq 0}  
n(e)\, 
t^{\vert n(e)\vert-1}=	\sum_{e\in \arr(\alpha)}  n(e)\, t^{\vert n(e)\vert-1}. $$
Substituting $t=1$, we obtain 
$$\sum_{k\geq 1} k\, u_k(\alpha)  =u'(1)=\sum_{e\in \arr(\alpha)}   n(e)=0.$$
The last equality follows from the fact that 	if an arrow $f$ links an 
arrow 
$e$ positively, then $e$ links $f$ negatively.	
	
		       \subsection{Examples}\label{sn:g23}  1.   For positive 
integers $p,q$, we define  $\alpha_{p,q}$ to be the lattice-looking 
virtual 
string formed by a   Euclidean circle in $\RR^2$ with counterclockwise 
orientation, $p$ 
disjoint vertical arrows $e_1,\ldots,e_p$ directed upward and numerated 
from 
left to right, and $q$ disjoint horizontal arrows   $e_{p+1},\ldots, 
e_{p+q}$  
crossing $e_1,\ldots,e_p$ from   right to left and numerated from 
bottom to top. 
(Here we identify arrows with geometric vectors in $\RR^2$  connecting 
two 
points of the core circle; the numeration of the arrows is compatible 
with the 
counterclockwise order of their tails.) Clearly, $n(e_i)= q$  for 
$i=1,\ldots , 
p$ and   
$n(e_{p+j})= -p$ for $j=1,\ldots , q$. Hence $u(\alpha_{p,q})=p t^{q}- 
q t^p$.  
We   conclude that the   strings $\{\alpha_{p,q}\}_{p\neq q}$ are 
pairwise 
non-homotopic and homotopically non-trivial.  The string    
$\alpha_{1,1}$ is 
homotopically trivial: it is obtained from a trivial  string by 
(b$)_s$. 
For $p\geq 2$, 
we have  $u(\alpha_{p,p})=0$. However 
$\alpha_{p,p}$ is homotopically non-trivial as will be shown  below. 

It follows from the definitions that $\overline {\alpha_{p,q}}=\alpha_{p,q}$ and  $ (\alpha_{p,q})^-=\alpha_{q,p}$.
Thus  the string $\alpha=\alpha_{p,q}$ with $p\neq q$ is not homotopic to  $\alpha^-, {\overline \alpha}^-$. 

Formula \ref{homr} implies that the strings  $\alpha_{p,1}$ and $\alpha_{1,p}$ with $p\geq 2$ have   minimal rank in
their homotopy classes.  We shall prove below that the same  holds for all  $\alpha_{p,q}$ with $(p,q)\neq (1,1)$.

2. A  permutation $\sigma$ of the set $\{1,2,\ldots , m\}$ gives rise to 
a 
virtual string $\alpha_\sigma$ of rank $m$ as follows. Let $S^1=\{z\in 
\CC\,\vert \,  \vert z\vert=1\}$ be the unit circle with 
counterclockwise 
orientation. For $i=1,\ldots ,m$, let $a_i $ (resp. $b_i $) be the 
point of 
$S^1$ with real part $(i-1)/m$ and  negative (resp. positive) imaginary 
part. 
Then $\alpha_\sigma$ is formed by $S^1$ and the  $m$ arrows
$\{(a_i, b_{\sigma (i)})\}_{i=1}^m$. For the $i$-th arrow  $e_i=(a_i, 
b_{\sigma 
(i)})$,
\begin{equation}\label{pup}n(e_i)=  \#\{j=i+1,\ldots, m\,\vert \,\sigma (j)< \sigma 
(i)\}-\#\{j=1,\ldots, 
i-1\,\vert \,\sigma (j)> \sigma (i)\}.\end{equation}
This allows us to compute the polynomial $u(\alpha_\sigma)$ directly 
from 
$\sigma$.  This example generalizes the previous one since $\alpha_{p,q}=\alpha_\sigma$ for the permutation  $\sigma$
of the set $\{1,2,\ldots , p+q\}$  given by
$$\sigma(i)=\left\{\begin{array}{ll}
i+q ,~ {\rm {if}} 
\,\,\, 
1\leq i\leq p \\
\noalign{\smallskip}
i-p
,~ 
{\rm
{if}} \,\,\, p<i\leq p+q.
\end{array} \right.$$

\subsection{Properties of $u$}\label{sn:g24} We point out a few simple 
properties of the polynomial $u$. For a virtual string $\alpha$, we 
have  
$u(\alpha)=u (\overline \alpha)$. This  follows from the fact that 
if two arrows are linked positively (resp. negatively), then 
the reversed arrows are also linked positively (resp. negatively).
The transformation $\alpha\mapsto   \alpha^-$  transforms positively 
linked 
pairs of arrows into negatively linked pairs and vice versa. Therefore
  $u (  \alpha^-)= - u (\alpha)$. 
As an 
application, 
we   observe that if $u(\alpha)\neq 0$, then $\alpha$ is  not homotopic 
  to 
$\alpha^-, {\overline \alpha}^-$.

It is obvious that if a  string $\alpha$ is a  product 
of 
 strings $\alpha_1$ and $\alpha_2$, then 
  $u (\alpha)= u (\alpha_1)+ u  (\alpha_2)$.

\begin{theor}\label{th:t21}
                     An integral polynomial $u(t)$ can be realized as 
the 
$u$-polynomial   of a virtual string if and only if $u(0)=u'(1)=0$.
                     \end{theor}
                     \begin{proof}
                     We need only to prove the sufficiency of the 
condition  
$u(0)=u'(1)=0$. The proof goes by induction on the degree of $u$. If 
this degree 
  is  $\leq 1$, then $u=0$ is realized by a trivial virtual string. 
Assume that 
our claim is true for   polynomials of   degree $<m$ where $m\geq 2$. 
Let $u(t)$ be 
a polynomial of degree $m $ with highest term $a t^m$ where $a\in \ZZ$ 
and 
$a\neq 0$. Then      $v(t)=u(t)-a (t^m-mt)$ is a polynomial of degree 
$<m$ with 
$v(0)=v'(1)=0$.	     By the inductive assumption, $v(t)$ is realizable as the $u$-polynomial   of a  string. 
By 
Example \ref{sn:g23}, the polynomial $t^m-mt$ is also realizable. Taking a product  of strings we observe that  the sum of
realizable polynomials is realizable.  Hence for 
$a>0$, the polynomial $u(t)=v(t)+a (t^m-mt)  $ is realizable. If $a<0$, 
then this argument shows that $-u(t)$ is realizable by a string, $\alpha$. Then
$u(t)$ is realized by $\alpha^-$.      \end{proof}
		     
\subsection{Computation for curves}\label{sn:g25} We compute the polynomial $u$ 
 for the   string
$\alpha=\alpha_\omega $  underlying a  closed curve $\omega:S\to \Sigma$ 
on a surface $\Sigma$. The computation goes in terms of  the homological 
intersection form $B:H_1(\Sigma)\times H_1(\Sigma) \to \ZZ$ determined 
by the 
orientation of $\Sigma$. Here and below $H_1(\Sigma)=H_1(\Sigma;\ZZ)$.

 Let $e=(a,b)$ be an arrow of $\alpha$.  
Then $\omega(a)=\omega(b)$ so that    $\omega$ transforms the 
arcs $  
ab, ba\subset S$   into   loops $\omega(ab), \omega(ba)$ in $\Sigma$. Set 
$[e]=[\omega(ab)]\in 
H_1(\Sigma)$ and $[e]^*= [\omega(ba)]\in H_1(\Sigma)$ where the square 
brackets on 
the right-hand side stand for the homology class of a loop. 
We compute the intersection number $B([e], [e]^*)\in \ZZ$. The  loops 
$\omega(ab), 
\omega(ba)$      intersect transversely 
except at their common origin $\omega(a)=\omega(b)$.  
Drawing a picture of   $\omega(ab), \omega(ba)$ in a neighborhood of    
$\omega(a)=\omega(b)$, 
one  observes that   a small deformation makes these   loops   
disjoint in 
this neighborhood. 
The  transversal intersections of $\omega(ab), \omega(ba)$ bijectively correspond to the arrows of 
$\alpha$ 
linked with $e$, i.e., the arrows connecting an interior point of  $ab 
$ with  
an interior point of   $ba $. The intersection sign at such an intersection  is $+1$    if the  tail of the corresponding  arrow lies
on $ab$ and is
$-1$  otherwise. 
Adding these signs, we obtain that  $B([e], [e]^*)=n(e) $. This formula can 
be 
rewritten in a more convenient form. Set   $   s=[\omega]= [\omega(S)]\in 
H_1(\Sigma)$. Observe 
that $ s=[e]+  [e]^*$ and therefore 
$$ B([e], [e]^*)=B([e], s-[e])=B([e],s )-
B([e],  [e])=B([e],s).$$
Thus
\begin{equation}\label{pik}n(e)=B([e],s).\end{equation}
Therefore  for any $k\geq 1$, $$u_k(\alpha)=\#\{e\in \arr(\alpha)\,\vert \,   B( [e],s)=k\}
- \#\{e\in \arr(\alpha)\,\vert \,   B([e],s)=-k\}  $$
and
$$ u(\alpha)=\sum_{e\in \arr(\alpha), B( [e],s)  \neq 0} \sign ( B( 
[e],s) )\, 
t^{\vert  B( [e],s)  \vert}.$$
Using the bijective correspondence between the set $\arr(\alpha)$  and the set $\Join\! \!  (\omega)$ of the double
points of   
$\omega$, we can rewrite the previous formula as  \begin{equation}\label{orrr} u(\alpha)=\sum_{x \in \Join  (\omega),
B( [\omega_x],s)  \neq 0}
\sign ( B(   [\omega_x],s) )\, 
t^{\vert  B(   [\omega_x],s)  \vert}\end{equation}
where  for $x \in \Join\! \! (\omega)$,  we let $\omega_x: [0,1]\to \Sigma$  be the loop beginning at  $x$
and following along $\omega$ until the first return to $x$
and such that the pair (a 
positive 
tangent vector of $\omega_x$ at $0$, a positive tangent vector of $\omega_x$ at 
$1$)  is a 
positive basis in   the tangent space of   $x$.  

\begin{theor}\label{th:t2547}
      The  polynomial  $u$ of a slice virtual string is equal to $0$.
                     \end{theor}
                     \begin{proof} Let $\alpha_0$ be  a slice virtual string  realized by a closed curve  $\omega_0$
on the boundary of a handlebody
$H$ such that $\omega_0$  is contractible in $H$.   By a {\it meridian} of $H$ we  mean an embedding
$S^1\hookrightarrow
\partial H$ which extends to an embedding of a 2-disc    into $H$.  Pick a base point $\ast\in \partial H$ and let  $1\in S^1$
be the base point of  
$S^1$. The kernel of the inclusion homomorphism
$\pi_1(\partial H,\ast)\to
\pi_1(H,\ast)$ is  normally generated   by the homotopy classes of    meridians   of $H$.  Therefore $\omega_0$ is
homotopic to a loop 
$\prod_{i=1}^n {\rrr}_i m_i{\rrr}^{-1}_i $ where   $m_1,...,m_n:S^1\hookrightarrow
\partial H$ are  meridians of $H$ and 
${\rrr}_i:[0,1]\to  \partial H$ is a path   leading from  $\ast$ to $m_i(1)$.  Deforming slightly these meridians and
paths we can assume that the images of 
$m_1,...,m_n$ are disjoint  simple closed  curves, that  ${\rrr}_1,...,{\rrr}_n$ meet these curves and each other
transversely, and that in a neighborhood of   $\ast$ the paths ${\rrr}_1,...,{\rrr}_n$ look like radii  going out of
$\ast$ in the cyclic order 
${\rrr}_1,...,{\rrr}_n$. Pushing   ${\rrr}_i$ slightly to the left (resp.  to the right) we obtain a \lq\lq parallel" path ${\rrr}_i^+$ (resp.
${\rrr}_i^-$). Doing it carefully we can assume  that  ${\rrr}_i^+ (0)= {\rrr}_{i-1}^-(0)$ for $i=1,...,n$ so that in a neighborhood of
$\ast$ the paths ${\rrr}_1^+, {\rrr}_1^-,..., {\rrr}_n^+,{\rrr}_n^-$ form $n$ disjoint embedded arcs approximating the $n$ radii above
from both sides.   We also assume that ${\rrr}_i^+(1)$ (resp. ${\rrr}_i^-(1)$) is a point of $m_i(S^1)$ lying just after (resp.
just before) $m_i(1)$.  Let $m'_i:[0,1] \to \partial H$ be an arc leading from  ${\rrr}_i^+(1)$ to  ${\rrr}_i^-(1)$ 
along  $m_i(S^1)$.  Then $ \omega={\rrr}_1^+ m'_1 ({\rrr}^{-}_1)^{-1}\cdots  {\rrr}_n^+ m'_n({\rrr}^{-}_n)^{-1}$  is a
generic loop in $\partial H$ homotopic to
$\omega_0$.  Let $\alpha$ be the underlying virtual string of $\omega$. The loop $\omega$ has two types of  
self-crossings:   each crossing of
${\rrr}_i$ with
$m_j$ (possibly $i=j$) gives rise to two self-crossings  $x, y$ of  $\omega$;      each
crossing of
${\rrr}_i$ with
${\rrr}_j$ gives rise to four self-crossings  $x,y,z,t$ of $\omega$.  We claim that in the first case $x,y$ contribute
opposite terms to  the right-hand side of  Formula \ref{orrr};  in the second case  the points $x, y,z,t$
can be partitioned into two pairs each contributing opposite terms to  the right-hand side of    \ref{orrr}.  This would imply
that $u(\alpha)=0$ and since $\alpha$ is homotopic to $\alpha_0$, we also have $u(\alpha_0)=0$.  To prove our claim in
the first case it suffices to check that      $B(  [\omega_{x}],s)=-B(  [\omega_{y}],s)$ where
$s=[\omega]\in H_1(\partial H)$.    It follows from the definitions of    $\omega_{x}, \omega_{y}$ that  
$ [\omega_{x}]+ [\omega_{y}]  =s\pm [m_i]$ for $i\neq j$ and  
$ [\omega_{x}]+ [\omega_{y}]  =s\pm (s-[m_i])$ for $i= j$.  The sign $\pm$ in these formulas depends on the intersection
sign at the       crossing of
${\rrr}_i$ with
$m_j$. Note that $B(s,s)=B([m_i],s)=0$ since $s=[m_1]+\cdots + [m_n]$   and the meridians   $m_1,...,m_n$ are disjoint. 
Hence  $B(  [\omega_{x}],s) +B(  [\omega_{y}],s)=0$. In the second case we can numerate the points $x,y,z,t$
so that  $ 
[\omega_{x}]+ [\omega_{y}]  =s+ [m_i]$,  $ 
[\omega_{z}]+ [\omega_{t}]  =s- [m_i]$ and   apply the same  argument as above.
\end{proof}

\begin{corol}\label{th:t254711}
      For any $p\neq q$, the string $\alpha_{p,q}$ is not slice.
                     \end{corol}

 \subsection{Remarks.}\label{sn:g26}  1. It is easy to verify that all 
virtual 
strings of rank $3$ are either homotopically trivial or  homeomorphic 
to     
$\alpha_{1,2},  \alpha_{2,1}$. Note   that the homotopy classes of  
$\alpha_{1,2}, \alpha_{2,1}$ are distinguished already by $u_1$; indeed $u_1(\alpha_{1,2})=-2$ and 
 $u_1(\alpha_{2,1})=2$. 

2.  For an   open string $\mu$ with core manifold $[0,1]$,  we can define two 
 polynomials
$u^+(\mu)$ and $u^-(\mu)$.    Observe that the set $\arr(\mu)$ of arrows of $\mu$ is a
disjoint union   $\arr^+(\mu)\cup \arr^-(\mu)$ where $\arr^{+} (\mu)$ (resp.  $\arr^{-} (\mu)$)  is the set of arrows
$ (a,b)\in \arr(\mu)$ with $a,b\in [0,1]$ such that $a<b$ (resp. $b<a$).   For $e \in \arr(\mu)$,   set $n(e)=n(e^{cl})\in
\ZZ$
 where  $e^{cl}$ is the corresponding arrow of the closure, $\mu^{cl}$, of
$\mu$. 
 For   $k\geq 1$,  set  $$u_k^{\pm} (\mu)=\#\{e\in 
\arr^{\pm}(\mu)\,\vert \,n(e) =k\}
- \#\{e\in \arr^{\mp}(\mu)\,\vert \,n(e) =- k\} \in \ZZ.$$
This number and the polynomials  $u^{\pm} (\mu)=\sum_{k\geq 1} u^{\pm}_k(\mu)\, t^k$ are homotopy invariants of
$\mu$.   Clearly, $u(\mu^{cl})=u^{+} (\mu)+u^{-} (\mu)$.  Using   $u^{\pm}$,  it is easy to give examples
of non-homotopic open strings with homotopic closures. 

3. For a virtual string $\alpha$ and a positive integer $p$, we define a 
virtual 
string $\alpha(p)$ as follows. Let us identify the core circle of 
$\alpha$ with 
$S^1\subset \CC$. Each arrow of $\alpha$ can be graphically  presented 
by a 
vector
in $\CC$ connecting two points of $S^1$. These vectors are mutually 
transverse. 
Now, replace each of these vectors, say $e$, by $p$   disjoint parallel 
vectors 
$e_1, \ldots , e_p$ running   closely to $e$ and having endpoints on  $S^1$. 
This 
gives  a virtual string $\alpha(p)$ of rank $pm$ where $m$ is the rank 
of 
$\alpha$. It is obvious that $n(e_1)=\ldots =n(e_p)= p \,n(e)$.
Therefore $u(\alpha(p)) (t)=p\, u (\alpha) (t^p)$. In particular, if 
$u(\alpha)\neq 0$, then $\alpha(p)$ is  
homotopically non-trivial. 

 \section{Geometric realization  of virtual strings}
                    
		         \subsection{Realization   of   strings}\label{fi:g31}  
We explain here that every virtual string   admits a canonical 
realization by a closed curve on a surface and moreover describe   all  its 
realizations.   

Let  $\alpha$ be a virtual string of rank $m$ with core circle $S$.    Identifying the  tail with the  head for all arrows 
of 
$\alpha$, we transform $S$ into   a 1-dimensional CW-complex  
$\Gamma=\Gamma_\alpha$.  We   thicken $\Gamma$ to a surface $\Sigma_\alpha$ as follows.
 If $m=0$, then 
$\Gamma=S$ and we set $\Sigma_\alpha= S\times [-1,1]$. Assume 
that $m\geq 1$.  The 0-cells (vertices) of $\Gamma$ have 
valency 4 
and their number is equal to  $m$. A neighborhood of a 
vertex 
$v\in \Gamma$ embeds into the unit 2-disc   $D^2=\{(x,y)\in 
\RR^2\,\vert \,  
x^2+y^2\leq 1\}$ as follows. Suppose that $v$ is obtained from an arrow 
$(a,b)$ 
where $a,b\in S$.  Note that any 
point 
$x\in S$ splits its small neighborhood in $S$ into two oriented arcs,  one 
of them being
incoming and the other one being outgoing with respect to $x$. 
Therefore a  neighborhood  of 
$v$ in 
$\Gamma$ consists of four  arcs which can be identified with small  
incoming 
and outgoing arcs of $a, b$ on $S$. We embed this neighborhood into $D^2$ so 
that $v$ 
goes to the origin and the incoming (resp. outgoing) arcs of $a,b$ 
go to 
the intervals $[-1,0]\times 0$,  $0  \times [-1,0]$ (resp. $[0,1]\times 0$, 
$0  \times 
[0,1]$), respectively. In this way all   vertices of 
$\Gamma$ are thickened 
 to copies of $D^2$ endowed with counterclockwise 
orientation. 
Each  1-cell  of $\Gamma$ connects two (possibly coinciding) vertices 
and is 
thickened  to a ribbon connecting the corresponding 2-discs. The 
thickening is 
uniquely determined by the condition that the orientation of these   
2-discs 
  extends to their union with the ribbon.  Thickening in this way 
all  the 
vertices and 1-cells of $\Gamma$ we embed $\Gamma$ into a 
surface  $ 
\Sigma_\alpha$.  By   
construction,  $\Sigma_\alpha$ is a   compact connected   oriented 
surface with non-void 
boundary and 
Euler characteristic   $\chi(\Sigma_\alpha)=\chi (\Gamma)=-m$. Composing the natural projection $S\to \Gamma$
with the  inclusion $\Gamma \hookrightarrow \Sigma_\alpha$ we obtain a 
closed 
curve   $\omega_\alpha:S\to \Sigma_\alpha$ realizing  $\alpha$.  
 The construction of $\Sigma_\alpha$ is well known, see 
\cite{f}, 
\cite{ca}, \cite{cw}.

			 It is clear that for any surface $\Sigma$ and any 
(generic) closed curve $\omega:S\to \Sigma$   realizing  $\alpha$, a  
regular  
neighborhood   of $\omega(S)$
in $\Sigma$ is homeomorphic to $\Sigma_\alpha$.  Moreover, the 
homeomorphism   can  be chosen to transform $\omega$ into $\omega_\alpha$. In other 
words, $\omega$ 
can be obtained as a composition of $\omega_\alpha$ with an 
orientation-preserving 
embedding $\Sigma_\alpha\hookrightarrow \Sigma$. In particular, $\Sigma_\alpha$ is 
a   
surface of minimal genus containing a closed curve realizing  $\alpha$.  Therefore the genus   $g(\alpha)$ of $\alpha$
defined in Section 
\ref{sn:g133} is equal to the genus of $\Sigma_\alpha$. It  will be  explicitly 
computed in 
the next subsection. Note finally that a  
closed 
surface of minimal genus containing   a curve   realizing  $\alpha$  is 
obtained 
from $\Sigma_\alpha$ by gluing 2-discs to all components of $\partial 
\Sigma_\alpha$. 
			 
			  \subsection{Homological computations}\label{fi:g32} 
Consider  again a virtual string $\alpha$ of rank $m$ with core circle $S$. 
Let   
$\Gamma=\Gamma_\alpha$,  $\Sigma=\Sigma_\alpha$, and $\omega=\omega_\alpha:S\to \Sigma_\alpha$ be the
graph,  the  surface,  and the closed curve    constructed in  the 
previous subsection. 	    The 
orientation of $\Sigma$ determines a homological intersection pairing   
$B=B_\alpha: 
H_1(\Sigma)\times H_1(\Sigma) \to \ZZ$. This bilinear pairing   is  
skew-symmetric  and its rank is equal to twice the genus of 
$\Sigma$. Thus
		 $$g(\alpha) = (1/2)\, {\text{rank}}  B_\alpha.$$ In particular, 
  $\alpha$ can be realized by a closed curve on $S^2$ or $\RR^2$  if and only if $B_\alpha=0$.

Since $\Gamma$ is a deformation retract of $\Sigma 
$,  the inclusion homomorphism $H_1(\Gamma)\to 
H_1(\Sigma 
)$ is an isomorphism. Since   $\Gamma$ is  a connected graph with  $\chi 
(\Gamma)=-m$, the 
group $H_1(\Gamma)=H_1(\Sigma 
)$ is a free abelian group of rank ${m+1}$. We  
describe a  canonical 
basis in $H_1(\Sigma 
)$. 
 Set $s= [\omega ] \in H_1(\Sigma)$.  For an arrow $e=(a,b)\in 
\arr(\alpha)$, the map $\omega$ transforms the arc $  ab\subset S$, leading 
from $a$ 
to $b$ in the positive direction, into a loop $\omega(ab)$ in $\Sigma$. Set 
$[e]=[\omega(ab)]\in H_1(\Sigma)$. An easy induction on $m$ shows  that 
$s\cup 
\{[e]\}_{e\in \arr(\alpha)} $ is a basis of $H_1(\Sigma)$. 
 Our next aim is to compute    the matrix of $B$ in this 
basis.   Note 
for the 
record that $B(x,y)=-B(y,x)$ and $B(x,x)= 0$ for all elements $x,y$ of 
this 
basis.

	By Formula \ref{pik},  $B([e],s)=n(e)$ for any $e\in 
\arr(\alpha)$. To compute  the other values of $B$, 
 we need more notation. Let $a,b$ be distinct point of  
$S$.   The {\it interior} of the arc  $ab\subset S$ is the set 
$(ab)^\circ= 
ab-\{a,b\}$. For  any   arcs $ab, cd \subset S$,  we define   $ab\cdot 
cd\in 
\ZZ$ to be the number of arrows of $\alpha$ with tail in 	
$(ab)^\circ$ and head in $(cd)^\circ$	    minus the number of arrows 
of 
$\alpha$ with tail in 	$(cd)^\circ$ and head in  
$(ab)^\circ$. 
Note that the arrows with both endpoints in $(ab)^\circ\cap (cd)^\circ$ 
appear 
in this expression twice with opposite signs and therefore cancel out. 
Clearly, 
$ab\cdot cd=-cd\cdot ab$.
 In particular,  $ ab \cdot   ab =0$. If   $e=(a,b)$ is 
an arrow 
of $\alpha$, then it follows from the  definitions that $n(e)=ab\cdot  ba$. 
More 
generally, for any arrow $f=(c,d) $ of $\alpha$ unlinked with 
 $e=(a,b)$, 
 \begin{equation}\label{equ01}n(e)=ab\cdot cd +ab \cdot dc.\end{equation}
 Note that the arrows $e=(a,b)$,  $f=(c,d)$  never contribute  to 
$ab\cdot cd$  and to   $ab \cdot dc$
because  neither   endpoint of $e$ lies in  $ (ab)^\circ$ and neither   
endpoint of 
$f$ lies in  $ (cd)^\circ$ or in $ (dc)^\circ$. 
Applying Formula \ref{equ01} to the string obtained from $\alpha$ by reversing the arrow $e$, we obtain that
 \begin{equation}\label{equ02}n(e) =-ba \cdot cd 
-ba \cdot 
dc.\end{equation}

		\begin{lemma}\label{l:t33} Let $e=(a,b)$ and $f=(c,d)$  be two 
  arrows of $\alpha$. Then $B([e],[f])=ab \cdot cd 
+\varepsilon$ where 
$\varepsilon=0$ if $e$ and $f$ are unlinked, $\varepsilon=1$ if   $f$ 
links 
$e$ positively, and  $\varepsilon=-1$ if   $f$ links $e$  negatively.
			   \end{lemma} 
                     \begin{proof}   If $e=f$, then $a=c, b=d$ and all terms of the stated equality are equal to 0.  (Note
that an arrow is unlinked  with itself.) Assume from now on that $e\neq f$ so that $a,b,c,d$ are  pairwise distinct points of
$S$.

Suppose first that
$e$ and
$f$ are unlinked. There are four cases to consider depending on whether the 
endpoints  of 
$e,f$ lie on $S$ in the cyclic order
(i)   $a,b, c,d$, or  (ii)   $a,b, d,c$, or  (iii)   $a,c,d, b$, or 
(iv)  
$a,d,c, b$.

		      In the case (i), the arcs $ ab , cd \subset S$ are 
disjoint so that   $[e],[f]\in H_1(\Sigma)$ are represented by transversal loops
$\omega(ab), 
\omega(cd)$, respectively. Then $B([e],[f])=ab \cdot cd $, cf. Section \ref{sn:g25}.
		     
		       In the case (ii), the arcs $ 
ab , dc \subset S$ are disjoint so that   $[e],[f]^{\ast}=s-[f]\in H_1(\Sigma)$ are 
represented by 
transversal loops $\omega(ab), \omega(d c)$, respectively. Hence $B([e],[f]^{\ast})= ab \cdot 
dc $  and $$B([e],[f])=B([e],s 
-[f]^{\ast})=B([e],s) - B([e],  [f]^{\ast}) =n(e) - ab \cdot 
dc= 
ab\cdot cd.$$
		     
		     In the case (iii), we have $B([e],[f])=-B([f],[e])=- c d 
\cdot ab=ab \cdot c d$ since the pair $(f,e)$ satisfies the conditions 
of (ii).
		     
		     In the case (iv),   the arcs $ 
ba , dc \subset S$ are disjoint so that   $[e]^{\ast}=s-[e], [f]^{\ast}=s-[f]\in H_1(\Sigma)$ are 
represented by 
transversal loops $\omega(ba), \omega(d c)$, respectively. 
Therefore 
$$B([e],[f])=B([e],s)+ B(s,[f])+B(s-[e],s-[f])=n(e)-n(f)+ B([e]^{\ast}, 
[f]^{\ast})
		     = n(e)- n(f) + ba \cdot dc.$$  
  It remains to observe that
		     $$ ba \cdot dc=-n(e) - ba \cdot  cd =-n(e)+cd\cdot ba= -n(e)+n(f) -cd\cdot ab=-n(e) +n(f) + ab \cdot 
cd.$$

		  Suppose that $f$ links $e$ positively. Then  their endpoints   
lie on $S$ in the cyclic order $a,c,b, d$.  The   loops 
$X=\omega(ab), Y=\omega(cd)$ representing $[e], [f]\in H_1(\Sigma)$ are not transversal since   both contain $\omega(c
b)$. Pushing 
$Y$  slightly  to its left in $\Sigma$, we obtain a loop, $Y^+$,  transversal to $X$. It is understood that the 
point 
$\omega(c)=\omega(d)\in Y$ is pushed  to a point lying between   $\omega(a c)$ and 
$\omega(d a)$	
in a small neighborhood of $\omega(c)=\omega(d)$.
Introducing coordinates $(x,y)$ in this neighborhood we can locally 
identify 
$X,Y, Y^+$ with the axis $y=0$, the union of two half-lines  $x=0, 
y\leq 0$ and 
$y=0, x\geq 0$, and the union
of two half-lines $x=-1, y\leq 1$ and $y=1, x\geq -1$, respectively. 
To compute the intersection number $B([e], [f])= X\cdot Y=X\cdot Y^+$, we  split the set $X\cap Y^+$ into 
four 
disjoint subsets. The first of them consists of a single point near 
$\omega(c)=\omega(d)$, given in the coordinates   above by $x=-1,y=0$. This point contributes $1$ to  $X\cdot Y^+$.
The second subset of $X\cap Y^+$ is   $\omega(ac)\cap Y^+$;  its points are numerated by  arrows of 
$\alpha$ with one endpoint  in the interior of $ac$ and the other 
endpoint in 
the interior of $cd$. The contribution of these crossings to 
$X\cdot Y^+$ is equal to $a c \cdot c  d$.
The third subset of $X\cap Y^+$ is numerated by   the   crossings of   
$\omega(cb)$ 
with the part of $Y^+$ obtained by pushing $\omega(bd)\subset Y$ to the 
left; they 
are numerated by arrows of $\alpha$ with one endpoint  in the interior 
of $c b$ 
and the other endpoint in the interior of $bd$. The contribution of 
these 
crossings to 
$X\cdot Y^+$ is   $c b \cdot bd$. The forth subset of $X\cap Y^+$ is 
numerated 
by   the   self-crossings of   $\omega(cb)$: each of them gives rise to two 
points of 
$X\cap Y^+$ with opposite intersection  signs. Therefore this forth subset 
contributes 0 to 
$X\cdot Y^+$.  Summing up these contributions we obtain
$$B([e],[f])=  1+  a c \cdot c  d + c b \cdot bd+0=a c \cdot c  d + 
cb\cdot cb + c 
b \cdot bd+1$$
$$=a c \cdot c  d + cb\cdot cd+1= a b\cdot cd+1.$$ 
	     
	If $f$ links  $e$ negatively, then $e$ links $f$ positively and by the 
results above,   $$B([e],[f])=-B([f],[e])=
		     -(c d \cdot a  b   +1)
		     =a  b \cdot c d -1.$$ \end{proof}			 
			    
     \subsection{Examples}\label{fi:g33}  (1) Consider the  string   
$\alpha=\alpha_{p,q}$ 
with $p,q\geq 1$ introduced in Section \ref{sn:g23}.1.  Recall the  arrows 
$e_1,\ldots, e_{p+q}$ of $\alpha$.  We  compute the  matrix of the
bilinear form $B=B_{\alpha}:H_1(\Sigma_\alpha)\times  
H_1(\Sigma_\alpha)\to \ZZ$ with respect to the basis $s\cup 
\{[e_i]\}_{i=1}^{p+q} $.  By Formula \ref{pik}, 
     $B([e_i],s) =q$ for $i=1,\ldots, p$ and $B([e_{p+j}],s)= -p$ for 
$j=1,\ldots , q$. Each pair of arrows $e_i, e_{i'}$  with $i, 
i'=1,\ldots, p$  is unlinked and by   Lemma \ref{l:t33}, $B([e_i],
[e_{i'}]) =0$. 
Similarly, each pair of arrows   $e_{p+j}, e_{p+j'}$ with   $j, 
j'=1,\ldots , q$ 
is unlinked and $B([e_{p+j}], [ e_{p+j'}])=0$. The arrow $e_{p+j}$ links 
$e_i$ 
positively and   by   Lemma \ref{l:t33},  $B([e_i],
[e_{p+j}])=(p-i)+(q-j)+1$. 
  It is easy to compute that the rank of $B$  is  equal to  2 
if 
$p=q=1$, to  6 if $\min(p,q)\geq 3$, and to $4$ in all the other cases.
 The genus $g(\alpha
) $, as we know, is half of this rank. 
 In particular, $g(\alpha_{1,1})=1$ which shows that the genus is not a 
homotopy 
invariant. 
 
 (2) Consider the  string   
$\alpha=\alpha_{\sigma}$ 
  defined  in Section \ref{sn:g23}.2  for a permutation $\sigma$ of the set $\{1,2,\ldots , m\}$. 
Recall the  arrows 
$e_1,\ldots, e_{m}$ of $\alpha$. 
We  compute the  
matrix of the bilinear form $B=B_{\alpha}:H_1(\Sigma_\alpha)\times  
H_1(\Sigma_\alpha)\to \ZZ$ with respect to the basis $s\cup 
\{[e_i]\}_{i=1}^m $. The number  $B([e_i],s) = n(e_i)$ is computed by Formula     \ref{pup}.
Pick two indices $i,j$ with  $1\leq i < j\leq m$.  Lemma \ref{l:t33} implies that if $\sigma(i)<\sigma(j)$, then
$$B([e_i], [e_j])=\#\{k\,\vert\, i<k<j, \sigma(j)<\sigma(k)\}\,-\, \#\{k\,\vert\, j<k\leq m,
\sigma(i)<\sigma(k)<\sigma(j)\}.$$
If $\sigma(j)<\sigma(i)$, then
$$B([e_i], [e_j])=\#\{k\,\vert\, i<k<j, \sigma(j)<\sigma(k)\}\,+\, \#\{k\,\vert\, j<k\leq m,
\sigma(j)<\sigma(k)<\sigma(i)\}+1.$$

 \subsection{Homotopy of strings re-examined}\label{sn:nb} The fact that all strings can be realized by curves on surfaces
allows us to reformulate the notion of homotopy of strings entirely in terms of homotopies of curves. For strings $\alpha,
\beta$, we write $\alpha\sim \beta$ if  these two strings can be realized by homotopic closed curves on the same surface.
The relation $\sim$ is reflexive and symmetric  but not transitive.  The next lemma shows that the relation of homotopy is 
precisely the equivalence relation generated by $\sim$. 

	\begin{lemma}\label{l:nbl} Two strings $\alpha,
\beta$ are homotopic if and only if there is a sequence of strings $\alpha_1=\alpha, \alpha_2,..., \alpha_n=\beta$ such that
$\alpha_i\sim \alpha_{i+1}$ for $i=1,..., n-1$.
			   \end{lemma} 
                     \begin{proof} As we know, the underlying virtual strings of homotopic closed curves on 
a surface are themselves homotopic.  Therefore if there is a sequence of strings $\alpha_1=\alpha, \alpha_2,...,
\alpha_n=\beta$ such that
$\alpha_i\sim \alpha_{i+1}$ for $i=1,..., n-1$, then $\alpha$ is homotopic to $\beta$.  To prove the converse, it suffices
to show that if $\beta$ is obtained from $\alpha$ by  a homotopy move  (a$)_s$, 
(b$)_s$,  or
(c$)_s$, then  $\alpha\sim \beta$.  

Let $S$ be the core circle  of $\alpha$ and $\omega:S\to 
\Sigma $ be 
a   curve    realizing  $\alpha$ on a surface $\Sigma$.  Pick   distinct points $a,b 
\in S$  
such that the arc   $ab\subset S$ does not contain  endpoints of 
$\alpha$.  
Let 
   $\beta$  be   obtained from $\alpha$ by the move (a$)_s$ adding    the arrow  $ (a,b)$.  
   Attaching to   $\omega$ a small curl 
on the 
right  of the arc $\omega(ab)$, we obtain a closed curve 
$\omega':S\to \Sigma $ realizing $\beta$.   Clearly, $\omega'$ is homotopic to $\omega$ in $\Sigma$.  Hence
 $\alpha\sim \beta$.

Pick two    arcs  $x,y$ on $S$ disjoint from each other and  
from the   
endpoints of $\alpha$. Let $a, a'$ be the endpoints of $x$ (in an 
arbitrary 
order) and $b,b'$ be the endpoints of $y$. Let $\beta$ be obtained from $\alpha$ by the move  (b$)_s$ adding 
to $\alpha$ the  arrows $ (a,b)$ and $ (b',a')$.    Let $D_x, D_y\subset
\Sigma - \omega(S)$ be two small closed discs  lying near 
the arcs $\omega(x), \omega(y)$, respectively. Removing the interiors of 
these discs 
from $\Sigma$ and gluing the circles $\partial D_x,\partial D_y$ along an 
orientation-reversing homeomorphism, we obtain a new (oriented) surface, 
$\Sigma'$, 
containing $\omega(S)$. In $\Sigma'$ the arcs $\omega(x)$ and $\omega(y)$ are 
adjacent 
to the component of 
 $\Sigma'-\omega(S)$ containing $\partial D_x=\partial D_y$.  We can push 
$\omega(x)$ 
across this component towards $\omega(y)$  and eventually across $\omega(y)$. This gives a  curve $\omega':S\to
\Sigma'$   realizing $\beta$ and homotopic to $\omega:S\to \Sigma'$. 
Hence
 $\alpha\sim \beta$. Note that the four 
possible 
forms of the   move   (b$)_s$ (depending on whether $x$ leads from $a$ to $a'$ or  from $a'$ to $a$   and similarly  for
$y$) are realized by choosing  
$D_x, D_y$  on the 
left or on the right of $\omega(x), \omega(y)$.

Suppose that $\alpha$ has three arrows 
$(a^+,b), (b^+,c), (c^+,a)$ where  $a, a^+, b,b^+$, $c, c^+\in S$ such that
the  (positively oriented) arcs   $aa^+$,   $bb^+$,    $cc^+$ are disjoint from each 
other and 
from the other   endpoints of $\alpha$. Let $\beta$ be obtained from $\alpha$ by the move  (c$)_s$ replacing  the 
arrows 
$(a^+,b), (b^+,c), (c^+,a)$ with the arrows  $ (a,b^+), (b,c^+), (c,a^+)$.
  Consider the canonical realization $ \omega_\alpha:S
\to
\Sigma_\alpha$ of $\alpha$.  Observe that  the    arcs $ \omega_\alpha  (aa^+),  \omega_\alpha (bb^+), \omega_\alpha  (cc^+)$     form  a simple
closed curve in
$\Sigma_\alpha$  isotopic to a boundary component of   $\Sigma_\alpha$. Gluing  a 2-disc  $D$
  to this boundary component we embed  $\Sigma_\alpha$ into a bigger 
surface, $\Sigma$. Pushing the branch   $ \omega_\alpha(aa^+)$ 
across $D\subset \Sigma$ 
and then across the double point $\omega(b^+)=\omega (c)$, we obtain a  curve $\omega':S\to
\Sigma$   realizing $\beta$ and homotopic to $ \omega_\alpha:S\to \Sigma$. 
Hence
 $\alpha\sim \beta$. 
\end{proof} 
  
			  \subsection{Adams operations}\label{sn:nbccc} We can define \lq\lq Adams operations" $\{\psi^n\}_{n\in \ZZ}$ on 
the  set $\str$  of homotopy classes of virtual strings.  Let $\alpha$ be a   string. Replacing $\alpha$ by a homeomorphic
string, we can identify its core circle with   
$S^1=\{z\in 
\CC\,\vert \,  \vert z\vert=1\}$. Consider  a   curve
$\omega:S^1\to
\Sigma$ realizing $\alpha$ on a surface $\Sigma$.  The   mapping  $  S^1\to \Sigma$ sending 
$z\in S^1$  to $ \omega(z^n)$    is homotopic to a generic curve $   S^1\to \Sigma$. We 
define 
$\psi^n(\alpha)\in \str$ to be the homotopy class of its underlying string. Lemma \ref{l:nbl} implies that   
$\psi^n:\str\to \str
$ is a well defined mapping. Clearly,
$\psi^{1} (\alpha)=\alpha$, $\psi^{-1} (\alpha)=\alpha^{-}$, and  $\psi^{mn}=\psi^m\circ 
\psi^n$ for any  $m,n\in \ZZ$.  As an exercise, the reader may check that $u(\psi^n(\alpha))=\sign (n)\, n^2\, u(\alpha)$.

			\section{Combinatorics of closed curves on the 2-sphere}\label{fi:g40}
                    
		         \subsection{Gauss words}\label{fi:g41} 
In this section, an {\it alphabet}  is a finite set and   
{\it 
letters} are its elements.  A {\it word} in an alphabet   is a finite 
sequence 
of letters. The words  will be  always considered up to 
circular 
permutations.  A {\it Gauss word}     is a word in an alphabet   in which all letters of the alphabet 
occur  
exactly twice.  Each virtual 
string gives rise to a Gauss word as follows.  Label the arrows of the  string     with 
different 
letters. Traverse  the core circle of the string in the  positive direction  and write down  the   label of  an arrow each time   
we  cross  its 
endpoint. This gives a Gauss word   well defined   up to circular 
permutations 
(and the choice of letters).  It is clear that the    
Gauss 
word associated   in the Introduction  with a  closed curve on $\RR^2$ (or more generally on any surface)    coincides with
the  Gauss 
word associated with  the underlying virtual string.

Gauss \cite{ga} studied the following question: when   
a Gauss 
word can be realized by a   closed curve on the plane~? He gave a 
necessary 
condition (condition (i) in Theorem \ref{th:e41} below) which he knew 
not to be 
sufficient. Later this question was studied by several authors, see 
\cite{cw} 
and M. L. Marx' review MR2000i:05056 for  references. A most    elegant 
solution to 
this question was given by P. Rosenstiehl \cite{ro}, see also \cite{rr} and   Corollary \ref{th:e415} below. 

Note   that a    Gauss  word      is realizable by a closed curve on 
$\RR^2$ 
if and only if it is realizable by a closed curve on the 2-sphere $S^2$. 
We shall focus on   curves on $S^2$ 
rather than 
on $\RR^2$.

		\subsection{Bipartitions}\label{fi:g4152}  A {\it 
bipartition} (or a {\it  partition into two sets})   of a set $E$ is 
a 
non-ordered pair of disjoint (possibly empty) subsets of $E$ whose union 
is $E$.     We explain now that   virtual
strings   and closed curves on surfaces naturally  give   rise   to  bipartitions. 
 Consider a 
virtual string   $\alpha$ with core circle $S$. For   arrows $e=(a,b)$ and $f=(c,d)$ of $\alpha$, we define  
$q(e,f)\in \ZZ $ to be   the    number of arrowheads  of
(the arrows  of) $\alpha$ lying on the semi-open arc  $ac-\{a\}\subset S$ minus  the   number of arrowtails  of
  $\alpha$ lying on   $ac-\{a\}$.     If   $e=f$, then by definition $q(e,f)=0$. 
 It is   easy to check that  $$q(e,f)+q(f,e)=0,\,\,\,\,\,\, q(e,f) +q(f,g)+q(g,e)=0$$ for any
    $e,f,g\in  \arr(\alpha)$.  We   use    $q\, (\modu 2)$ to define  an equivalence
relation on
$\arr(\alpha)$: two arrows 
$e,f\in \arr(\alpha)$ are  {\it  equivalent}  if $q(e,f)\equiv 0\, (\modu 2)$. 
This can be reformulated in simpler terms:  the arrows $e=(a,b)$ and $f=(c,d)$  are  
equivalent if either $e=f$ or the number   of    
endpoints 
of $\alpha$ lying in the interior of  the arc $ac\subset S$ is odd. 
It is obvious that this   relation on
$\arr(\alpha)$  has at most 
two 
equivalence classes. They form a   partition of $ \arr(\alpha)$ into 
two subsets 
 (one of  them may be empty).		
If the   arrows  of $\alpha$ are labeled by  different letters, then 
this 
bipartition of $\arr(\alpha)$ induces a bipartition of the set of 
letters.  Thus $\alpha$ gives rise   to a pair (a Gauss word, a bipartition of the alphabet).

Applying this construction   to the  underlying virtual string  of a  closed curve $\omega$ on a surface, we obtain 
a bipartition  of the  set of  double points of $\omega$.		If the  double points of $\omega$  are labeled by  different
letters, then  this 
    induces a bipartition of the set of 
letters.  Thus $\omega$ gives rise   to a pair (a Gauss word, a bipartition of the alphabet).
For  curves on $S^2$, we shall give a  geometric   re-formulation of this bipartition in Remark \ref{fi:g44}.1.

	\subsection{Curves on $S^2$}\label{fi:g40zz}   It is natural to ask whether the pair  (the Gauss
word, the bipartition of the alphabet) associated with a curve on a surface  is a   full  homeomorphism invariant of the curve
and what    values it can take.  We shall answer these  questions for curves on
$S^2$.

	We begin with  definitions. Let $w$ be a   word in an alphabet $E$. 
We 
say that two  (distinct) letters $i,j\in E$ are {\it $w$-interlaced} if $w$ has the 
form 
$i...j...i...j...$   up to a circular  permutation. For   $i\in E$, denote 
by $w_i$ 
the set of letters $w$-interlaced with $i$.   A bipartition of $E$ is {\it compatible}  with $w$ if it satisfies the following
condition:  for any   $w$-interlaced $i,j\in E$, we have   $\#(w_i \cap 
w_j)\equiv 0 \, (\modu 2)$  if  $i,j$  belong to different subsets of the 
bipartition   and $\#(w_i \cap 
w_j)\equiv 1 \, (\modu 2)$  if   $i,j$  belong to the same subset.
			 
\begin{theor}\label{th:e41}  A pair (a  Gauss  word $w$ in an   
alphabet $E$, a 
bipartition of $E$)  is realizable by a   closed curve on $S^2$ if 
and only if 
the following three conditions are satisfied:

(i)  for all $i\in E$, the set $w_i$ has an even number of elements: 
$\#(w_i)\equiv 0 \, (\modu 2)$;

(ii) if $i,j\in E$ are not $w$-interlaced, then $\#(w_i \cap w_j)\equiv 0 \, (\modu 2)$;

(iii)  the bipartition of $E$ is compatible with $w$.
 \end{theor}
		     \begin{proof} We begin with   preliminary computations.  
Let $\alpha$ be a virtual string of rank $m\geq 1$ with core circle 
$S$.   
Suppose that the arrows of $\alpha$ are labeled by   elements of  a set 
$E$  so that each   $i\in E$ appears as a label of a unique arrow 
$e_i=(a_i,b_i)$ of $\alpha$  where $a_i,b_i\in S$. As explained above, 
this 
gives  a Gauss word $w=w(\alpha)$ in the alphabet $E$. 	
Observe that    letters 
$i,j\in E$ are $w$-interlaced if and only if $e_i,e_j$ are linked.  
Therefore  
for any $i\in E$, the number $\#(w_i)$ is the number of arrows of 
$\alpha$ 
linked with $e_i$. Hence $$\#(w_i)\equiv n(e_i) \, (\modu 2).$$ Similarly, for any $i,j\in E$, we have   
$\#(w_i \cap 
w_j) \equiv n_{i,j}\, (\modu 2)$ where $n_{i,j}$ is the number of  arrows of $\alpha$
  linked with both $e_i$ and $e_j$.  We now  relate $n_{i,j}$  to 
the 
number $B([e_i], [e_j]) $ defined in Section \ref{fi:g32}.  We  
shall assume   that $n(e_i)\equiv n(e_j)\equiv 0 \, (\modu 2)$ (this will be enough for our aims).  Applying Formulas
\ref{equ01},  
\ref{equ02}  to
$e=e_i, e_j$ we obtain that   \begin{equation}\label{equ1} a_i b_i\cdot a_j  b_j \equiv  a_ib_i \cdot b_ja_j \equiv b_ia_i
\cdot a_jb_j  
\equiv b_ia_i 
\cdot b_j a_j \,\,(\modu 2).\end{equation}
If $i,j $ are not 
$w$-interlaced, 
then   one of the arcs $a_ib_i, b_i a_i$, say  $\beta$, is disjoint 
from one of 
the arcs $a_jb_j, b_j a_j$, say $\gamma$.   It is obvious that  
	 $n_{i,j}\equiv \beta \cdot \gamma \, (\modu 2)$.  Lemma \ref{l:t33} and Formula  \ref{equ1} 
 imply that 
 \begin{equation}\label{equ3}\#(w_i \cap w_j)\equiv n_{i,j} \equiv \beta \cdot \gamma  \equiv  a_i b_i\cdot a_j  b_j
 =  B([e_i], [e_j]) \, (\modu 2).\end{equation} Suppose now  that  $i,j$ are $w$-interlaced. It will be enough for our
aims to consider the case where 
 the endpoints of $e_i,e_j$ lie on $S$ in the cyclic order 
$a_i,a_j,b_i, b_j$. 
Then  $n_{i,j}\equiv  a_i a_j   \cdot 
b_ib_j +  
a_j b_i \cdot b_j a_i \, (\modu 2)$.   
 Let  $r$     be the number of endpoints of 
$\alpha$ 
lying on the arc   $a_j b_i$.   It is easy 
to see 
that $r\equiv a_j b_i \cdot b_i a_j \, (\modu 2)$.  The number $q(e_i,e_j)\in \ZZ $ defined in  Section
 \ref{fi:g4152} satisfies 
$$  q(e_i,e_j)+r  \equiv \#(w_i)
\equiv n(e_i)\equiv 0 \, (\modu 2).$$ Therefore $  q(e_i,e_j)\equiv r \equiv  a_j b_i \cdot b_i a_j \, (\modu 2)$.
 Using Lemma \ref{l:t33}, we obtain  
  \begin{equation}\label{equ2}B([e_i],[e_j])  =a_ib_i \cdot a_jb_j +1 = 
a_ia_j 
\cdot a_jb_i +a_ia_j \cdot  b_i b_j+ a_j b_i \cdot b_i 
b_j+1 \end{equation}
$$\equiv  - a_ia_j 
\cdot a_jb_i +a_ia_j \cdot  b_i b_j+ a_j b_i \cdot b_i 
b_j+2 \,a_j b_i \cdot b_j a_i+ 1
$$
$$= (a_jb_i \cdot  a_ia_j  + a_j b_i \cdot b_i b_j+  a_j b_i \cdot b_j a_i ) + 
(a_i a_j  
 \cdot b_ib_j +  a_j b_i \cdot b_j a_i)+ 1$$
$$\equiv a_j b_i \cdot b_i a_j +n_{i,j} +1 \equiv  q(e_i,e_j)+\#(w_i \cap 
w_j)+1\, (\modu 2). $$

		      We   can now prove the necessity of  the conditions (i) -- (iii)  of the theorem.   Suppose that a Gauss word 
$w$ in an alphabet  $E$ is realized by a  closed curve   $\omega: S^1\to 
S^2$. Let 
$\alpha$ be the underlying virtual string of $\omega$. Thus, 
$w=w(\alpha)$ for 
an appropriate bijective labelling $E\to  \arr (\alpha), i\mapsto e_i$.  
   By Formula \ref{pik},   $n(e_i)$ is the   
intersection number 
of two cycles lying in a neighborhood of 
	 $\omega(S^1)$ in $S^2$.  Since the   intersection number of any  two 
cycles  in $S^2$ is zero,  $n(e_i)=0$ and  $\#(w_i)\equiv  n(e_i) \equiv 0 \, (\modu 2)$ for all $i\in E$. 
Condition (ii)  follows  similarly from Formula \ref{equ3}.
It remains to verify  that the bipartition of $E $ induced by $\alpha$ is compatible with $w$. Let 
$i,j\in E$ be  $w$-interlaced. Permuting if necessary $i$ 
and $j$ we 
can assume that 
 the endpoints of $e_i=(a_i,b_i)$, $e_j=(a_j,b_j)$   lie on $S^1$ in the cyclic order 
$a_i,a_j,b_i, b_j$. 
  Formula \ref{equ2} implies that  $ q(e_i,e_j)+\#(w_i \cap 
w_j)+1
\equiv 0\, (\modu 2)$.   
 Hence
 $\#(w_i \cap w_j)\equiv 0\, (\modu 2)$ if and only if $q(e_i,e_j)\equiv 1\, (\modu 2)$, i.e., 
if and 
only if $i,j$ belong to different equivalence classes.

To accomplish the proof of the theorem,  we need a  general construction  of strings from  Gauss
words and  bipartitions.  Assume that we have   a Gauss word
$w=z_1z_2\ldots z_{2m}$    in an alphabet
$E$ with 
$\#E=m$ 
and a bipartition  $E=X\cup Y$.  Assume that $w$ satisfies the \lq\lq 
parity 
condition" (i) of  the theorem. We construct  a string 
$\alpha $   giving rise to this word and this bipartition. 
 The core circle of $\alpha$ is  the circle 
$S=\RR\cup 
\{\infty\}$ with orientation extending the right-handed orientation on 
$\RR$. 
The set of arrow endpoints of $\alpha$ is   the set  $\{1,2,\ldots, 
2m\}\in 
\RR\subset S$. This set is partitioned into $m$   pairs: two   
points 
$a,b\in \{1,2,\ldots, 2m\}$ form a pair if $z_a=z_b$. We   order each 
such pair 
$\{a,b\}$ as follows. By the parity condition (i), we have $a-b\equiv 1\, (\modu 2)$. 
 Therefore one of the numbers $a,b$ is even and the other one is odd.  
If 
$z_a=z_b\in X$, then we put the odd one on the first place and the even 
one on 
the second place. If $z_a=z_b\in Y$, then we do the opposite. The  string $\alpha$ is the circle $S$ with these 
$m$ ordered pairs  of points. It follows from 
the 
definitions that the  Gauss word of $\alpha$ is $w$ and the induced 
bipartition of 
$E$ is $E=X\cup Y$. 
 
To establish the sufficiency of the conditions (i) -- (iii) of the theorem,  we need only  to show that if $w$ satisfies (i),  (ii)
and  the  bipartition $E=X\cup 
Y$   is compatible with $w$, then the string   $\alpha $ 
constructed 
in the previous paragraph  is realizable by a closed curve in $S^2$. It 
suffices to 
show that the surface $\Sigma=\Sigma_\alpha$ is a disc with holes; then 
$\Sigma$ 
embeds in $S^2$ so that the canonical realization of $\alpha$ in 
$\Sigma$ 
gives a realization in $S^2$. By the classification of compact  
surfaces, it suffices to 
prove that 
the intersection form $B: H_1(\Sigma)\times H_1(\Sigma) \to \ZZ$ takes 
only even 
values.
 By Formula \ref{pik} and Lemma   \ref{l:t33}, we need only to prove 
that the 
expressions for the values of $B$ appearing there   are   even. For all $i\in E$, 
we have   $n(e_i)\equiv \#(w_i)\equiv 0 \, (\modu 2)$ by  the condition  (i). Let  $i,j$ be distinct elements of
$E$.    If 
the  arrows $e_i, e_j\in \arr(\alpha)$ are unlinked, then  
  $B([e_i],[e_j]) \equiv \#(w_i \cap w_j)  \equiv 0 \, (\modu 2)$ 
by  Formula \ref{equ3} and the condition (ii). 
Suppose that $e_i,e_j$ are linked. 
Assume first that  the endpoints  of $e_i=(a_i,b_i), e_j=(a_j,b_j)$ lie on $S$ in the cyclic 
order 
$a_i,a_j,b_i, b_j$.
 If $i,j\in X$ (resp. $i,j\in Y$), then both numbers  $a_i,a_j \in  \{1,2,\ldots, 2m\}$ are   odd  (resp. even)
so that 
 $q(e_i,e_j)\equiv a_i-a_j \equiv 0\, (\modu 2)$. By   (iii), we have  $\#(w_i \cap w_j)\equiv 1 
\, (\modu 2)$. By Formula \ref{equ2}, $B([e_i],[e_j])$ is even.  If $i,j$ lie in different 
subsets of the bipartition, then 
 one of the numbers $a_i,a_j$ is even and the  other one is odd so that 
 $q(e_i,e_j)\equiv a_i-a_j \equiv 1\, (\modu 2)$. By    (iii),  $\#(w_i \cap w_j)\equiv 0 \, (\modu
2)$
 and   by Formula \ref{equ2}, $B([e_i],[e_j])$ is even. The case where the   endpoints  of $e_i, e_j$
lie  on $S$ in the cyclic order $a_i,b_j,b_i, a_j$ follows from the previous 
one 
using the skew-symmetry of   $B$. \end{proof}
 			 
\begin{corol}\label{th:e415} (P. Rosenstiehl  \cite{ro}) A Gauss  word $w$ in an  alphabet 
$E$ is 
realizable by a   closed curve in $S^2$ if and only if  $w$
satisfies the 
conditions 
(i),  (ii)  of Theorem \ref{th:e41} and   $E$ admits a bipartition 
compatible with $w$. \end{corol}

  We say 
that two  closed curves  $\omega,\omega':S^1\to \Sigma$  on a surface $\Sigma$ are {\it 
homeomorphic} 
if there is a homeomorphism $\varphi:\Sigma\to \Sigma$ such that 
$\omega'=\varphi\circ \omega$.  Here  $\varphi$ is not required to preserve 
orientation in $\Sigma$. 

     \begin{theor}\label{th:e431} Two closed curves  on 
$  S^2$   yield the same  pair (a  Gauss  word,       a   bipartition of the alphabet) if and only if these curves are
homeomorphic.
\end{theor}
 
                     \begin{proof}  We first prove that homeomorphic closed curves  on $S^2$ 
give rise 
to the same Gauss words and the same bipartitions. This is obvious if  the homeomorphism $S^2\to S^2$  relating the
curves is   orientation-preserving.  It remains to prove   that  the Gauss word and the bipartition  associated with a closed curve
$\omega:S^1\to S^2$ are preserved under reversal of  orientation of $S^2$.  By its very definition, the Gauss word does not
depend on this orientation.  It remains  to check that the string $\alpha$ underlying $\omega$ and the inverse string
$\overline
\alpha$ give rise to the same bipartition of the set  $\arr(\alpha)=\arr (\overline \alpha)$.  Let $e=(a,b), f=(c,d)\in
\arr(\alpha)$. An easy count shows  that the number of endpoints of $\alpha$ in the interior  of the arc
$bd\subset S^1$ is equal modulo 2  to the number of endpoints of $\alpha$  in the interior  of the arc $ac\subset S^1$
plus
$n(e)+n(f)$.  Since $\alpha$ underlies a curve on the 2-sphere,
$n(e)\equiv n(f)\equiv 0 (\modu 2)$. This shows that $e,f \in \arr(\alpha)$ are equivalent  in the sense of Section
\ref{fi:g4152}   if and only if  the  corresponding arrows
 $ (b,a), 
 (d,c)\in
\arr(\overline \alpha)$ are equivalent. 
  
Conversely, suppose that  two closed curves 
$\omega,\omega':S^1\to S^2$   yield the same  pair (a  Gauss  word in 
an alphabet $E$,       a   bipartition of $E$). We shall show that $\omega,\omega'$ are homeomorphic.   Let 
$\alpha , \alpha'$ be the underlying virtual strings of   $\omega,\omega' $, respectively.   Our assumption implies
that we can label the arrows of  $\alpha, \alpha'$ with  elements of  $E$ so that $\alpha, \alpha'$  give rise to the
same  Gauss word 
and the same bipartition of $E$. Since $\alpha, \alpha'$  determine the 
same Gauss word,   we may assume that they
 coincide up to  the choice of 
orientation   of  arrows.  We claim that either $\alpha=\alpha'$ or  these strings are opposite on all
arrows.  Suppose that $\alpha$ has an arrow
$  (a,b)$ such  that 
$(b,a)\in \arr(\alpha')$.  If $(c,d)$ is another arrow of $\alpha$, then  as in the previous paragraph, a simple
count  shows that   the number of endpoints of $\alpha$ lying inside the arc $ac$ has the 
opposite    parity to 
the number of endpoints of $\alpha$ lying inside the arc $bc$.  Since $\alpha, \alpha'$  determine the 
same 
bipartition of the set of arrows,  the pair $(c,d)$ cannot be an arrow of $\alpha'$.  Hence $(d,c)\in \arr(\alpha')$
which proves the claim above.  
 (This claim is compatible with the fact that the set of orientations on $m$
arrows has
$2^m$  elements while   the set of bipartitions of $E$ has $2^{m-1}$ elements where  $m=\#(E)$.)   Thus 
$\alpha'$ is homeomorphic  either to $\alpha$  or to its inverse   
$\overline{\alpha}$.   The second case can be reduced to the first 
one by 
composing $\omega$ with  an orientation-reversing 
self-homeomorphism of $S^2$. Thus we can assume that $\alpha$, 
$\alpha'$ are 
homeomorphic.    As we know,  $\omega$ is  a composition of   
the canonical realization $\omega_{\alpha}:S^1\to 
\Sigma_{\alpha}$ with an embedding  $\Sigma_{\alpha}\hookrightarrow S^2$ and similarly for $\omega'$. A 
homeomorphism $\alpha\to \alpha'$ extends to a homeomorphism $\varphi  : 
\Sigma_{\alpha}\to \Sigma_{\alpha'}$ transforming $\omega_{\alpha}$ 
into 
$\omega_{\alpha'}$. Note that all the components of $S^2-\Sigma_{\alpha} , S^2-\Sigma_{\alpha'}$ 
are 
discs. Since any homeomorphism of circles extends to a 
homeomorphism 
of discs bounded by these circles, $\varphi $ extends to a homeomorphism 
$S^2\to 
S^2$  transforming $\omega$ into $\omega'$. \end{proof}
		   
Theorem  \ref{th:e431}  shows that  the Gauss word and the bipartition associated with a closed curve on
$S^2$ is a full homeomorphism invariant.
  Theorems \ref{th:e41} and \ref{th:e431} give a complete  combinatorial description of  the set of homeomorphism
classes of curves on
$S^2$ in
 terms of Gauss words and bipartitions. 

  \subsection{Examples}\label{fi:g42} The   parity condition (i) of Theorem \ref{th:e41}  was 
pointed out 
by Gauss who knew that from $m=5$ on it is not sufficient. He gave as 
examples 
the sequences $1231245345$ and $1231435425$. They   satisfy (i)  but 
not 
(ii). The word  $w=123456214365$  satisfies (i) and (ii) but the 
alphabet
$\{1,2,3,4,5,6\}$ does not admit a bipartition compatible with $w$. 
  Indeed,  the letters  $1,3,5$ of this alphabet  are pairwise $w$-interlaced and   $ \#(w_1 \cap w_3)=\#(w_1\cap
w_5)=\#(w_3\cap w_5)=2 $. The compatibility would imply that $1,3,5$   belong to pairwise  different subsets 
of the 
bipartition which is impossible.

 \subsection{Irreducible Gauss words.}\label{fi:g43ff}     A Gauss word  is  {\it irreducible} if neither itself nor  its  circular
permutations can  be written as a concatenation of two  non-empty Gauss  words (in smaller alphabets). For example, the
word $1212$ in the  alphabet 
$\{1,2\}$ is irreducible while the word $1221$ is not. The next theorem shows that  if an    irreducible  Gauss words allows 
a
compatible bipartition of the alphabet   then such a bipartition  is unique.

      \begin{theor}\label{th:e431cdcd}  For an  irreducible Gauss word $w$ in an alpabet $E$, there is at most one bipartition
of
$E$ compatible with $w$.   \end{theor}
		     \begin{proof}    Consider the   graph $P$   whose vertices  are elements of $E$ 
and in 
which two vertices   are connected by an edge if and only if they are 
$w$-interlaced.      Since $w$ is irreducible,   $P$ is connected. 
(If it 
were disconnected, then realizing $w$ by a virtual string with   core 
circle  
$S^1\subset \CC$ and arrows represented by   geometric vectors, one 
would easily 
observe  that $w$ is not irreducible.) If a bipartition   $E=X\cup Y$ 
is compatible with $w$, then knowing for an element of $  E$ whether it lies in $X$ or 
$Y$ we 
can determine this   for   its immediate neighbors in   $P$. 
Proceeding   along $P$, we eventually determine this for all 
elements  of 
$E$. Hence there is at most one such bipartition of $E$.
  \end{proof}

	 \begin{corol}\label{th:e466} (C. H. Dowker and   M. B. Thistlethwaite \cite{dt})
Two closed curves   in $S^2$ realizing the same irreducible  
Gauss  word 
    are  homeomorphic. \end{corol}
 
This follows directly from Theorems \ref{th:e431} and \ref{th:e431cdcd}.

	 \begin{corol}\label{th:e466sd} 
If a 
Gauss word  
$w$  in an alphabet $E$ is obtained by concatenation of $k\geq 1$ irreducible 
Gauss 
words, then the number of  bipartitions  of $E $ compatible with $w$ is either 0 or  $2^{k-1}$. \end{corol}

\begin{proof}  Denote the set of
bipartitions of the alphabet compatible with  $w$ by $B(w)$.  Let
$w$ be obtained by concatenation of two Gauss words $w_1, w_2$ in disjoint alphabets $E_1,E_2$, respectively. 
Intersecting a bipartition of  $E=E_1\cup E_2$ with $E_1, E_2$ we obtain a  mapping, $j$,  from the set of bipartitions of
$E$ into the set of pairs (a bipartition of $E_1$, a bipartition of $E_2$).  This mapping   is 2-to-1. Since the
letters of
$E_1$ are not
$w$-interlaced with the letters of
$E_2$, we have  $j(B(w))=B(w_1) \times B(w_2)$. 
 Hence $\#(B(w))=2 \, \#(B(w_1)) \, \#(B(w_2))$.  This implies our claim by induction on $k$, the case  $k=1$ being 
Theorem
\ref{th:e431cdcd}. \end{proof}
  
    \subsection{Remarks}\label{fi:g44}  1.  If  a closed curve $\omega: S^1\to \Sigma$  on a surface
$\Sigma$  is $\ZZ/2\ZZ$-homologically trivial,  then  the  associated bipartition of the set of its double points $\Join\!  
(\omega)$ admits a simple geometric interpretation.   The  homological triviality of $\omega$ implies that  the
components of
$\Sigma-
\omega(S^1)$ can be colored white or black so that   the    components of $\Sigma-
\omega(S^1)$ adjacent to the same arc in
$\omega(S^1)-\Join\! (\omega)$  from opposite sides  have different colors.   With each   point $x\in \Join\! (\omega)$ we
associate    the component   of   $\Sigma-
\omega(S^1)$ adjacent to $x$ and  lying between the  two positive tangent vectors of
$\omega$ at
$x$.  It is easy to deduce from the definitions that two points of $\Join\! (\omega)$ belong to the same subset of the
bipartition determined by $\omega$  if and only if the associated  components of    $\Sigma-
\omega(S^1)$ have the same color.  In particular,  for such $\omega$, the  bipartition of
$\Join\! (\omega)$  determined by $\omega$  does not depend on the orientation of $\Sigma$.   

    2.  Theorem \ref{th:e431} can be generalized  to certain curves on surfaces of  arbitrary genus. Let us call a curve
 $\omega: S^1\to \Sigma$  on a surface $\Sigma$ {\it special} if  it is  $\ZZ/2\ZZ$-homologically trivial and all the
components of $\Sigma-
\omega(S^1)$ are discs.  The   proof of Theorem  \ref{th:e431}  shows that two special closed
curves  on  a closed surface   yield the same  pair (a  Gauss  word,       a   bipartition of the alphabet) if and only if these
curves are homeomorphic.

3.  To study curves on $\RR^2$   one can  involve an additional  piece of combinatorial data which is
a   subset of the alphabet. Indeed,  the complement of   a curve  in $\RR^2$ has one infinite region; the
labels of the  double points adjacent to this region form a subset  of the alphabet.    It is easy to see that the
Gauss word, the bipartition, and this  subset
form a full homeomorphism invariant of a closed curve on $\RR^2$.

  \section{Based  skew-symmetric  matrices}
                    
We introduce   algebraic notions  used in the sequel to define further  homotopy invariants  of   strings.

\subsection{Definitions}\label{fi:g51}   A  {\it based skew-symmetric matrix over 
$\ZZ$}  or shortly a {\it based matrix}  is a triple $(G, s, b)$ where   $G$ is a finite set, $s\in G$, 
and 
$b:G^2=G\times G \to \ZZ$ is a skew-symmetric mapping  so   that  
$b(g,h)=-b(h,g)$ for all $g,h\in G$. In particular,  $b(g,g)=0$ for all 
$g\in 
G$.    

Two based   matrices $(  G, s, b)$ and $(  G', s', b')$ are {\it 
isomorphic} 
if there is a bijection $G\to G'$  sending $s$ into $s'$  and  
transforming $b$ 
into $b'$.  To specify the isomorphism class of a   based matrix  $(G, s, b)$, it suffices  to specify the
matrix
$(b(g,h))_{g,h\in G}$ where  it is understood that the first  column and row correspond to $s$. In this way every
skew-symmetric square matrix over $\ZZ$ determines a  based matrix.

We   call an element $g\in G-\{s\}$ {\it annihilating} (with respect to 
$b$) if 
$b(g,h)=0$ for all $h\in G$. We call $g\in G-\{s\}$ a {\it core 
element} if 
$b(g,h)=b(s,h)$ for all $h\in G$.   We   call two elements $g_1,g_2\in 
G-\{s\}$ 
{\it complementary} if $  b (g_1, h)+ b
 ( g_2,h)=  b(s,h)$ for all $h\in  G$.  A based matrix 
$( G, s, b)$ is  {\it primitive}    if it has no annihilating elements, no core 
elements, and 
no complementary pairs of elements. An example of a primitive based matrix 
is 
provided by the {\it trivial based matrix} $(G,s, b)$ where $G$ consists of 
only one 
element $s$ and $b(s,s)=0$.

\subsection{Equivalence of based matrices}\label{fi:g515} We define 
three 
operations $M_1, M_2, M_3$ on based matrices, called {\it elementary 
extensions}. They add to a  based matrix $(  G, s, b)$ an 
annihilating 
element, a core element, and  a pair  of complementary elements, 
respectively. 
More precisely, $M_1$ transforms $(  G, s, b)$ into the (unique) based matrix $(  
\overline G=G\amalg \{g\},  s, \overline b)$
such that    $\overline b:\overline G\times \overline G \to \ZZ$  extends
$b$ and  $\overline b (g,h)= 0$ for all $h\in \overline G$.  
The move  $M_2$  transforms $(  G, s, b)$ into the (unique) based matrix $(  
\tilde 
G=G\amalg \{g\},   s, \tilde b)$
such that    $\tilde  b:\tilde  G\times \tilde  G \to \ZZ$  
extends $b$ and  $\tilde  b (g,h)=\tilde  b (s,h)$ for all $h\in 
\tilde G$. 
 The   move  $M_3$  transforms $(  G, s, b)$ into a   based matrix 
$(  \hat 
G=G\amalg \{g_1,g_2\},  s, \hat b)$
where   $\hat b:\hat G\times \hat G \to \ZZ$ is any skew-symmetric   map   
extending 
$b$ and such that $\hat b (g_1, h)+\hat b
 ( g_2,h)=   b(s,h)$ for all $h\in \hat G$.   It is clear that a based matrix 
$( G, s, b)$ is primitive  if  and only if it cannot be obtained from another 
based matrix by  an elementary extension.

Two based matrices are {\it   homologous}
if one can be obtained from   the other by a finite sequence of 
elementary 
extensions $M_1,M_2, M_3$, the inverse transformations, and 
isomorphisms.   The homology is an equivalence relation on the set of based matrices. A simple homology invariant of
a  based matrix
$T=(G,s,b)$ is provided by the 1-variable polynomial
  $$ u_T(t)=\sum_{e\in G, b(e,s)  \neq 0} \sign ( b(e,s) )\, 
t^{\vert  
b(e,s)  \vert}.$$

\begin{lemma}\label{l:t51} Every based matrix is  obtained from   a 
primitive based matrix by elementary extensions. Two homologous primitive based matrices are isomorphic.
		   \end{lemma}
                     \begin{proof}  The first claim is obvious: 
eliminating    
annihilating elements,   core elements, and   complementary pairs of 
elements by 
the moves $M_i^{-1}$ with $i=1,2,3$ we can  transform any 
based matrix $T$ into a  primitive based matrix $T_0$.  Then  $T$ is 
obtained from 
$T_0$ by elementary extensions. 
  
  To prove the second claim, we need the following assertion:
  
  $(\ast)$ a  move $M_i$ followed by   $M_j^{-1}$ yields the same 
result as    
an isomorphism, or a   move $ M_k^{\pm 1}$, or     a  move $M_k^{-1}$ 
followed 
by $M_l$ with $k,l\in \{1,2, 3\}$.

  This assertion will imply the second claim of the lemma. Indeed,    
suppose that two primitive  based matrices $T,T'$ are related by a finite 	
sequence 
of transformations $M_1^{\pm 1}, M_2^{\pm 1}, M_3^{\pm 1}$ and 
isomorphisms. 	
 An isomorphism of based matrices followed by $ M_i^{\pm 1}$     
can be also 
obtained as  $ M_i^{\pm 1}$    followed by an isomorphism.  Therefore 
all 
isomorphisms in our sequence   can be   accumulated at the end. The 
claim  
$(\ast)$ implies that $T,T'$
can be  related by a finite sequence of moves consisting of  several 
  moves of type $  M_i^{-1}$ followed by  several 
  moves of type $ M_i$ and   isomorphisms. However, since   $T$ is 
primitive  
  we cannot apply to it a move of type $  M_i^{-1}$. Hence there are no 
such 
moves in our sequence. Similarly, since $
  T'$ (and any isomorphic based matrix) is primitive, it cannot be obtained by 
an 
application of $M_i$. Therefore our sequence  consists solely of 
isomorphisms so 
that $T$ is isomorphic to $T'$.

 Let us now prove $(\ast)$. We have to consider  nine  cases depending on 
the values 
$i,j\in \{1,2,3\}$. 
 
  For $i, j \in \{1,2\}$, the move   $M_i$ on a  based matrix 
$(G,s,b)$ adds 
one element $g$   and then  $M_j^{-1}$ removes   one  element  $g'\in G\amalg \{g\}$. If 
$g'=g$, 
then   $  M_j^{-1}\circ M_i$   is the identity. If   $g'\neq g$, then 
$g'\in G$ 
is annihilating (resp. core) for $j=1$ (resp. $j=2$). The 
transformation $  
M_j^{-1}\circ M_i$   can be achieved by first applying $M_j^{-1}$ that 
removes 
$g'$ and then applying $M_i$ that adds $g$.

   Let $i=1, j=3$. The move 
 $M_i$ on $(G,s,b)$ adds an annihilating element $g$  and  
  $M_j^{-1}$ removes two complementary  elements $g_1,g_2\in G\amalg \{g\}$. 
  If $g_1\neq g$ and $g_2\neq g$, then $g_1,g_2\in G$  and $  
M_j^{-1}\circ M_i$ 
 can be achieved by first removing $g_1,g_2$ and then adding $g$. If 
$g_1=g$, 
then   $g_2 $ is a core element of $G$ and $  M_j^{-1}\circ M_i$ is the 
move 
$M_2^{-1}$    removing $g_2$. The case $g_2=g$ is similar.
  
   Let $i=2, j=3$. The move 
 $M_i$ on $(G,s,b)$ adds a core element $g$  and  
  $M_j^{-1}$ removes two complementary  elements $g_1,g_2\in G\amalg \{g\}$. 
  If $g_1\neq g$ and $g_2\neq g$, then  $  M_j^{-1}\circ M_i$  can be 
achieved 
by first removing $g_1,g_2$ and then adding $g$. If $g_1=g$, then   
$g_2\in G$ 
is   an annihilating element of $G$ and $  M_j^{-1}\circ M_i$ is the 
move 
$M_1^{-1}$    removing $g_2$. The case $g_2=g$ is similar.

   Let $i=3, j=1$. The move 
 $M_i$  on $(G,s,b)$ adds two complementary  elements $g_1,g_2$   and  
$M_j^{-1}$ removes an annihilating element  $g\in G\amalg \{g_1,g_2\}$. 
  If $g\neq g_1$ and $g\neq g_2$, then $g\in G$ and $  M_j^{-1}\circ 
M_i$  can 
be achieved by first removing $g$ and then adding $g_1,g_2$.
   If $g=g_1$, then   $ g_2$ is a core element of $G\amalg \{g_2\}$   
and $  
M_j^{-1}\circ M_i=M_2$. The case $g=g_2$ is similar.
  
   Let $i=3, j=2$. The move 
 $M_i$  on $(G,s,b)$ adds two complementary  elements $g_1,g_2$   and  
$M_j^{-1}$ removes a core element  $g\in G\amalg \{g_1,g_2\}$. 
  If $g\neq g_1$ and $g\neq g_2$, then $g\in G$ and  $  M_j^{-1}\circ 
M_i$  can 
be achieved by first removing $g$ and then adding $g_1,g_2$.
   If $g=g_1$, then   $ g_2$ is  an annihilating  element of $G\amalg 
\{g_2\}$   
and $  M_j^{-1}\circ M_i=M_1$. The case $g=g_2$ is similar.

      Let $i=j=3$. The move 
 $M_i$  on $(G,s,b)$ adds two complementary  elements $g_1,g_2$   and  
$M_j^{-1}$ removes two  complementary elements  $g'_1,g'_2\in G\amalg \{g_1,g_2\}$. If these 
two pairs  
 are disjoint, then 
  $  M_j^{-1}\circ M_i$  can be achieved by first removing  $g'_1,g'_2\in G$ 
and then 
adding $g_1,g_2$. If these two pairs  
  coincide, then $  M_j^{-1}\circ M_i$  is the identity. It remains to 
consider 
the case where these pairs have one common element, say $g'_1=g_1$, 
while 
$g'_2\neq g_2$. Then $g'_2\in G$ and for all $h\in G$,
  $$\hat b(g_2,h)=  b(s,h)-\hat b(g_1,h)=  b(s,h)-\hat b(g'_1,h)=\hat 
b(g'_2,h)= 
b(g'_2,h).$$
 Therefore the move   $  M_j^{-1}\circ M_i$ gives a based matrix isomorphic 
to 
$(G,s,b)$. The isomorphism $G\to (G-\{g'_2\}) \cup \{g_2\}$ is the 
identity on 
$G- \{g'_2\}$ and  sends $g'_2$ into $g_2$.
\end{proof}

Lemma \ref{l:t51}  implies that each  based matrix $T=(G,s,b)$ is homologous to a   primitive    based matrix $T_0=(G_0,s,b_0)$
unique up to  isomorphism.  This  reduces the  classification of homology classes of 
based matrices to an isomorphism classification of primitive based matrices. Note that   we can choose 
$T_0 $ in its isomorphism class so that $G_0\subset 
G$ and 
$b_0$ is the restriction of $b$ to $G_0\times G_0$.  

Each   isomorphism invariant $v$ of primitive  
based matrices  
 extends to a  homology        invariant   of based matrices by  $v(T)=v(T_0)$.   
The most important numerical  invariant of a primitive based matrix 
$(G,s,b)$ 
is the number  $\#(G)$. It is   easy to define further  
invariants  of primitive  based matrices. For instance, for   $k\in \ZZ$,  we can 
set
$$v_k(G,s,b)= \#\{g\in G\,\vert\, b(g,s )=k\}.$$
Similarly, for $k\in \ZZ$ and   a finite set of integers $A$ endowed 
with 
non-negative multiplicities,
set 
$$v_{k,A}(G,s,b)=\#\{g\in G\,\vert\,  b(g,s)=k\,\,{\text {and}}\,\,\, 
\{b(g,h)\}_{h\in G-\{s\}}=A\},$$
where the latter equality is understood as an equality of sets with 
multiplicities.  Clearly, $v_k=\sum_A v_{k,A}$.
 			 
\subsection{Remark}\label{fi:g52} The moves $M_1, M_2, M_3$ on
based matrices are 
not independent. It is easy to present $M_2$ as a composition of $M_3$ 
with 
$M_1^{-1}$.			 

\subsection{Genus of   based matrices}\label{fi:g5999} We define a numerical invariant of a  based matrix 
$T=(G,s,b)$ called its {\it genus} and denoted $\sigma (T)$.  For    subsets $X,Y\subset G$, set $b(X,Y)=\sum_{g\in X,
h\in Y} b(g,h)\in \ZZ$.  Clearly,  $b(X,Y)=-b(Y,X)$  and  $b(X,X)=b(\emptyset ,X)=0$ for all $X,Y\subset G$. A  {\it
regular partition} of
$G$ is  a splitting of
$G$ as a union of disjoint (possibly empty) subsets
$\{X_i\}_i$ such that $\# (X_i) \leq 2$ for all $i$ and one of $X_i$ is  the one-element set $\{s\}$.  The {\it matrix} of the
regular partition $x=\{X_i\}_i$   is the  matrix $(b(X_i,X_j))_{i,j}$.  This is a skew-symmetric square matrix over $\ZZ$. Its
rank is   an even number; let  $\sigma(x)$ denote half of this rank.  By definition, $\sigma(T)=\min_x
\sigma(x) $ where $ x$ runs over all regular partitions of $G$.   Extending $b$ by linearity to the lattice
generated by $G$ and identifying subsets of $G$ with vectors in this lattice whoose coordinates are   0 or 1, we can
interpret
$\sigma(T)$ as half the minimal rank of the restriction of $b$  to the sublattices arising from   regular partitions of $G$.

Note that $\sigma (T)\geq 0$ and  $\sigma (T)=0$ if and only if $G$ has a regular partition $\{X_i\}_i$ such that
$b(X_i,X_j)=0$ for all $i,j$.  In the latter case we say that   $T$ is {\it hyperbolic}.  It is easy to see that if  $T$ is
hyperbolic, then $u_T=0$.

The key property of the genus $\sigma (T)$ is
contained in the following lemma. 

\begin{lemma}\label{l:gg1} The genus of a  based matrix is  a  homology  invariant.
		   \end{lemma}
                     \begin{proof} 
  By Remark \ref{fi:g52}, it suffices to prove that $\sigma(T)=\sigma(T')$ for any  based matrix
$T'=(G',s,b')$     obtained from a  based matrix 
$T=(G,s,b)$ by a move
$M_i$ with $i=1,3$.  The set $  G'-G$ consists of one  element if $i=1$ and of two elements if $i=3$.  Pick a   
regular partition $x=\{X_i\}_i$ of
$G$ such that  $\sigma(T)= 
\sigma(x)$.  Consider the   regular partition $x'=(G'-G)\cup  \{X_i\}_i$ of $G'$.  Its matrix   is obtained from the one of
$x$ by adjoining a row and a column.  For $i=1$, these row and column are zero so that   $  \sigma
(x)
=
\sigma(x')$.   For $i= 3$,  we have $b(G'-G,Y)=b(\{s\}, Y)$ for all  $Y\subset G$.  Since one of the sets $X_i$ equals
$\{s\}$,  we again obtain  $  \sigma
(x)
=
\sigma(x')$.   Hence
$ \sigma(T')\leq \sigma(x') =\sigma(x)=\sigma(T) $. 

To prove the opposite inequality, pick  a   
regular partition $x'=\{X_i\}_i$ of
$G'$ such that  $\sigma(T')= 
\sigma(x')$.  Consider  first the case $i=1$.  One of the sets $X_i$ contains the 1-element set $G'-G$.  We replace this
$X_i$ by $X_i -(G'-G)$ and keep all the other $X_i$. This gives a regular partition  $x$   of $G$ whose matrix coincides
with the matrix  of $x'$.  Hence  $\sigma(T)\leq 
\sigma(x) =  
\sigma(x')= \sigma(T')$.  Let now $i=3$.  If  one of the sets $X_i$ is equal to $G'-G=\{g_1,g_2\}$, then  removing this
$X_i$ from
$x'$ we obtain a regular partition  $x$  of $G$.  As in the previous paragraph, 
$\sigma(x)=\sigma(x')$.   Hence  $\sigma(T)\leq 
\sigma(x) =  
\sigma(x')= \sigma(T')$.  Suppose that the   elements $g_1,g_2$  of $G'-G$ belong to different subsets, say $X_1,
X_2$,  of the partition
$x'$.   Then the sets $X_i$ with $i\neq 1,2$ and $X= (X_1\cup X_2) - \{g_1,g_2\}$ form a  regular partition of $G$. 
Let $X_0$ be the  term of the partitions $x$ and  $x'$ equal to $\{s\}$.  For
any $Y\subset G$, we have 
$$b(X,Y)=b(X_1,Y)+b( X_2,Y)- b(g_1,Y)- b(g_2,Y)=b(X_1,Y)+b( X_2,Y)- b(X_0,Y).$$ 
   Therefore the skew-symmetric  bilinear form  determined by the  matrix
of $x$ is induced from  the  skew-symmetric   bilinear form  determined by the  matrix
of $x'$ via the  linear map of the corresponding lattices sending the vectors $X$ and $\{X_i\}_{i\neq 1,2} $ respectively to 
$ X_1+X_2-X_0$ and $\{X_i\}_{i\neq 1,2}$.  Hence  
$\sigma(x)
\leq   
\sigma(x')$ and    $\sigma(T)\leq 
\sigma(x) \leq   
\sigma(x')= \sigma(T')$.
\end{proof}

\begin{corol}\label{aaal:gg1} A   based matrix   homologous to a hyperbolic based matrix is itself hyperbolic.
		   \end{corol}

\subsection{Transformations  $T\mapsto -T$, $ T\mapsto T^-$}\label{fi:g51253}   
 We define  two more operations on based   matrices.   For  a  based matrix $T=(G,s,b)$, set  $-T= (G,s,-b)$ and 
$T^-=(G,s, b^-)$  where     $b^- (s,h)=-b( s,h)$, $b^- (h,s)=-b( h,s) $  for all $h\in G$ and
$b^-(g,h)= b(g,h)+ b(s,g)-b(s,h) $ for all
$g,h\in G-\{s\}$.   The transformations  $T\mapsto -T$, $ T\mapsto T^-$ are commuting involutions on the set
of based matrices.   It is easy to  check that they are compatible
with   homology and preserve the class of   primitive based matrices.
It follows from the definitions that   $(-T)_0 =- T_0$ and $(T^-)_0=(T_0)^-$.  

\subsection{Remark}\label{uu:g52} We can define the  direct  sum  $T_1\oplus T_2$ of based matrices
$T_1=(  G_1, s_1, b_1)$, 
$T_2=(  G_2, 
s_2, b_2)$ to be   the based matrix $(  G, s, b)$   where  $G=(G_1\amalg  G_2 )/s_1=s_2$,  the element $s\in G$ is
defined by $ s=s_1=s_2$, and  
$b:G^2 \to \ZZ$ extends  both $b_1$ and $b_2$   and satisfies  $b(g_1,g_2)=0$ for any $g_1\in 
G_1-\{s_1\}$,  $g_2\in G_2-\{s_2\}$.  As an exercise, the reader may check that  the  direct  sum of primitive based matrices is
primitive and  the based matrix
$T\oplus (-T)$ is hyperbolic for any
$T$. 

  \section{Based matrices  of  strings}

\subsection{The based matrix of a string}\label{fi:g53} With each 
virtual 
string $\alpha$ we associate a  based matrix 
$T(\alpha)=(G,s,b)$. Set $G=G(\alpha)=\{s\}\amalg 
\arr(\alpha)$. To 
define $b=b(\alpha):G\times G\to \ZZ$, we identify $G$ with  the basis $s\cup 
\{[e]\}_{e\in \arr(\alpha)} $ of 
$ H_1(\Sigma_\alpha)$, see 
Section \ref{fi:g32}.  The map $b$ is obtained by restricting the 
homological 
intersection pairing
$ H_1(\Sigma_\alpha)\times H_1(\Sigma_\alpha)\to \ZZ$
to $G$. It is clear that $b$ is skew-symmetric. We can   compute 
$b$  combinatorially using   Formula \ref{pik} and Lemma  
\ref{l:t33}. In 
particular,    $b(e,s)=n(e)$ for all $e\in \arr(\alpha)$. 

The map  $b$ can be computed from any closed curve $\omega $ realizing $\alpha$ on a 
surface $\Sigma$. 
Indeed,   such a  curve is  obtained from  the 
canonical realization of $\alpha$ in $\Sigma_\alpha$  
via 
 an orientation-preserving embedding  $\Sigma_\alpha\hookrightarrow \Sigma$. 
It remains   to observe that such an embedding preserves  
intersection numbers    and  transforms the basis $s\cup 
\{[e]\}_{e\in \arr(\alpha)} $ of 
$ H_1(\Sigma_\alpha)$ into the subset $[\omega  ], \{[\omega_x] \}_{x\in \Join (\omega)}$ of $ H_1(\Sigma)$, cf.
Section
\ref{sn:g25}.

 \begin{lemma}\label{th:e53}  If two virtual strings   are homotopic, 
then their 
based matrices   are homologous. \end{lemma}
  \begin{proof} By Lemma \ref{l:nbl} it is enough to show that if  two closed curves $\omega, \omega' $ on a
surface
$\Sigma$ are homotopic, then the based matrices of their underlying strings are homologous.  
By the discussion in Section \ref{sn:g13},  it suffices to consider the  case where $\omega'$ is obtained from
$\omega$ by one of the   local moves listed there. 		

  If $\omega'$ is obtained from $\omega$ by   
adding a small  curl, then $\Join\! (\omega')=\Join\! (\omega)\cup \{y\}$  where $y$ is a new crossing. 
Clearly
$[\omega'_y]=0\in H_1(\Sigma)$  or  $[\omega'_y]=[ \omega']=[\omega] \in H_1(\Sigma)$ depending on
whether  the curl lies on the  right or  on the left  of
$\omega$. Also      $[\omega'_x]=[\omega_x]$  for all $x\in \Join\!(\omega)$. 
 Hence   
$T(\beta)$ is obtained from  
$T(\alpha)$ by $M_1$ or $M_2$.

Suppose that  $\omega'$  is obtained from $\omega$ by  the move  pushing a branch of $\omega$ across  another branch
and creating  two new double points $y,z$.  Clearly,   $[\omega'_x]=[\omega_x]$  for all $x\in \Join\!(\omega)\subset
\Join\!(\omega')$.  It is easy to see that  
 $[\omega'_y]+[\omega'_z]= [ \omega']=[\omega]\in H_1(\Sigma)$.
Therefore 
 $T(\beta)$ is obtained from  $T(\alpha)$ by $M_3$.
 
If  $\omega'$  is obtained from $\omega$ by   pushing  a branch of $\omega$ across a 
double point, then the
  subsets $[\omega  ], \{[\omega_x]_{x\in \Join (\omega)} \}  $ and $[\omega'  ],
\{[\omega'_x]_{x\in \Join (\omega')} \}  $ of $ 
H_1(\Sigma)$ coincide so that   $T(\alpha)$ is isomorphic to $T(\beta)$.  \end{proof} 
  
    \subsection{Homotopy invariants of strings from 
based matrices}\label{fi:g54}   
Every virtual string $\alpha$ gives rise 
to a   
primitive based matrix $T_0(\alpha) $  by $T_0(\alpha)=(T(\alpha))_0$. This is the only primitive
based matrix (up to  isomorphism) homologous to
$T(\alpha)$.   By
 Lemma  \ref{th:e53}, the based matrix  
$T_0(\alpha)=(G_0,s,b_0)$   is a 
  homotopy invariant of    $\alpha$.  This based matrix  determines the polynomial $u(\alpha)$ introduced in Section
\ref{abric}:  it follows from    Formulas \ref{homr11} and  \ref{pik} that 
  $  u(\alpha)=u_{T(\alpha)}=u_{T_0(\alpha)} $.  The number
$\rho (\alpha) =\#  (G_0 ) -1$ 
is a useful  homotopy invariant of $\alpha$
  which may be non-zero even when $u(\alpha)=0$, cf.  the 
examples 
below.  Note that  if    $\alpha$ is homotopically trivial, 
then 
$T_0(\alpha)$ is a trivial based matrix and $\rho(\alpha)=0$.

 It follows from the definitions that $T(\alpha^-)= (T(\alpha))^-$ and therefore  $T_0(\alpha^-)= (T_0(\alpha))^-$.
Similarly, $T(\overline \alpha)= -(T(\alpha))^-$ and    $T_0(\overline \alpha)= -(T_0(\alpha))^-$.

   The based matrix  $T_0(\alpha)=(G_0,s,b_0)$ can be used to estimate the 
homotopy 
rank and the  homotopy genus of $\alpha$. Namely,
   $hr(\alpha)\geq \rho(\alpha)$ since any string homotopic to $\alpha$ 
must 
have at least $\rho(\alpha)$ arrows. Similarly,  $hg(\alpha)\geq (1/2) 
\rank 
b_0$ where   $\rank 
b_0$  is     the 
rank of the integral matrix $(b_0(g,h))_{g,h\in G_0}$.   Indeed, if $\alpha'$ is a string homotopic to $\alpha$ and $T 
(\alpha')=(G',s',b')$, then  
   $g(\alpha')= (1/2) \rank b' \geq (1/2) \rank b_0$ since the matrix 
of $b'$ 
contains the matrix of $b_0$ as a submatrix. 

Combining the  inequalities   $hr(\alpha)\geq \rho(\alpha)$, 
$hg(\alpha)\geq (1/2) 
\rank 
b_0$  with the obvious inequalities $\rank  \alpha \geq
hr(\alpha)$ and
$g(\alpha)\geq hg(\alpha)$, we obtain that  if  
$T(\alpha)  $ is  primitive, then  $hr(\alpha) =\rank  \alpha$ and 
$hg(\alpha)  =g(\alpha)$.
  
 The next theorem gives an estimate for the  slice genus $sg(\alpha)$ of $\alpha$ via 
$\sigma(T(\alpha))=\sigma(T_0(\alpha))$.

 \begin{theor}\label{th:vv53}  For any string $\alpha$, we have $\sigma (T(\alpha))\leq 2\, sg(\alpha)$. \end{theor}
  \begin{proof}  Set
$k=sg(\alpha)$.  We  can  present
$\alpha$ by a loop on the boundary of a    handlebody $H$  which bounds a  (singular)  surface  of genus $k$  in  $H$. This
loop
  is homotopic in $\partial H$  to a loop  
$$\omega= \prod_{i=1}^n {\rrr}_i^+ m'_i ({\rrr}^{-}_i)^{-1} \prod_{j=1}^k p^+_j q^+_j (p^-_j)^{-1} (q^-_j)^{-1} $$
where ${\rrr}_i^+, {\rrr}^{-}_i,  m_i,  m'_i $ are  as in  the proof of Theorem \ref{th:t2547}  and the paths $p^-_j,  q^-_j  
$ on $  \partial H$ 
  are obtained from paths $p^+_j , q^+_j   $ on $ \partial H$, respectively,   by slight  pushing to the  right.  Choosing  
the  paths $p^{\pm}_j,  q^{\pm}_j, r^{\pm}_i$  carefully, we  can assume that  they begin and end in a small  disc 
$V\subset  \partial H - \cup_{i=1}^n m_i
$    and have no crossings in $V$.  Then the  crossings of 
$\omega$   split into pairs of points $(y,z)$ arising when 

 (a) the paths
$p^+_j, p^-_j$ meet  one of the  paths $  q^\pm_{k} $; 

 (b) the paths
$p^+_j, p^-_j$ meet  one of the  paths $p^\pm_{k} , {\rrr}^\pm_i, m'_i$; 

(c) the paths
$q^+_j, q^-_j$ meet  one of the  paths $q^\pm_{k}, {\rrr}^\pm_i, m'_i$; 

(d) the paths
$r^+_i, r^-_i$ meet  one of the  paths ${\rrr}^\pm_l$;

(e) the paths
$r^+_i, r^-_i$ meet  one of the  paths $ m'_l$. 

Such  pairs $(y,z)$ give rise to homology classes  $[\omega_{y}], [ \omega_{z}] \in
H_1(\partial H)$ whose sum can be explicitly computed.  Set   
 $s=[\omega ]=[m_1]+...+[m_n]\in
H_1(\partial H)$ and let $[p_j], [q_j]\in H_1(\partial H, V)=
H_1(\partial H)$ be the homology classes of  the paths  $p^\pm_j, q^\pm_j$, respectively.  
In the case (a),   $[\omega_{y}]+[ \omega_{z}]=s\pm [q_j]$ if $k\neq j$ and   $[\omega_{y}]+[ \omega_{z}]=s\pm (s-
[q_j])$ if $k=j$.  In the case (b),  $[\omega_{y}]+[ \omega_{z}]= s\pm
[q_j]$.   In the case (c),  $[\omega_{y}]+[ \omega_{z}]= s\pm
[p_j]$.   In the case (d),    $[\omega_{y}]+[ \omega_{z}]= s\pm
[m_i]$.    In the case (e),      $[\omega_{y}]+[ \omega_{z}] =s\pm [m_i]$ if $l\neq i$ and    $[\omega_{y}]+[
\omega_{z}]= s\pm (s- [m_i])$ if $l=i$.     These computations show that  the  sublattice of  $H_1(\partial H)$  generated
by such sums 
$[\omega_{y}] + [\omega_{z}]$ is contained in the sublattice   of  $H_1(\partial H)$   generated by  $2k+n$ elements  
$[p_j], [q_j],  [m_i]$. Since
$ B([m_i], [m_{l}])=0$ for all
$i,l$,  the restriction of  the intersection form  $B:H_1(\partial H)\times  H_1(\partial H)\to \ZZ$ to the latter  (and hence to
the former)  sublattice   has rank
$\leq 4k$.  Therefore for the underlying string  $\alpha_\omega$  of $\omega$ we have  $  \sigma(T(\alpha_\omega))\leq
2k
$. Since 
$\alpha_\omega$ is homotopic to $\alpha$,  their based matrices are homologous. By Lemma \ref{l:gg1},   
$ 
\sigma(T(\alpha))=\sigma(T(\alpha_\omega)) \leq 2k$. 
\end{proof}

\begin{corol}\label{th:vwww53}  For a  slice string $\alpha$, the based matrices  $ T(\alpha) $ and $T_0(\alpha)$ are hyperbolic.
\end{corol}
    
     \subsection{Applications}\label{fi:g55}  (1) The based matrix 
$T(\alpha_{p,q})$ 
of the string   $\alpha_{p,q}$ with $p,q\geq 1$ was computed in Section 
\ref{fi:g33}. It is easy to check   that except in the case $p=q=1$, this 
based matrix  
is primitive. Thus $ T_0(\alpha_{p,q})=T(\alpha_{p,q})$,  $hr(\alpha_{p,q})=\rank \alpha_{p,q}=p+q$  
and  $hg(\alpha_{p,q})=
 g(\alpha_{p,q})$ provided $p\neq 1$ or
$q\neq 1$.  In particular,   $\alpha_{p,p}$
  is a homotopically non-trivial  string with zero $u$-polynomial 
for all 
$p>1$. 
 
 (2)  We  show  that the product of strings defined in Section \ref{sn:g24} does not
induce  a well-defined operation on the set of homotopy classes of strings.  To this end, we exhibit a homotopically
trivial string  
  whose product with itself is not homotopically trivial.  Consider the permutation $\sigma=(12)(34)$ on the set
$\{1,2,3,4\}$ permuting 1 with 2 and 3 with 4.  Consider the rank  4  string
$ \alpha_\sigma$, as defined in Section \ref{sn:g23}.2.  Drawing a picture of $\alpha_\sigma$,
one  observes that  it  is a product of two copies of 
$ \alpha_{(12)}$.  The latter string   is   homotopically trivial since it is obtained from a trivial string by the
homotopy move (b$)_s$.  The based matrix  $T({\alpha_\sigma})$
  can be explicitly computed, cf. Section    \ref{fi:g33}.2.   It is determined by the
following  skew-symmetric   matrix:  
$$
 \left [ \begin{array}{ccccc}  0& -1&1& -1& 1 \\
        1& 0& 1& -1& 1 \\
-1& -1& 0& -1&1 \\
1& 1&1& 0&1 \\
-1& -1& -1& -1& 0 
\end{array} \right ].$$
It is easy to check that this based matrix is primitive. Hence $\alpha_\sigma$ is not homotopically trivial.  Moreover,  it is not
homotopic to a string with   $<4$ arrows.

(3) We   prove that the involution $\alpha\mapsto \overline \alpha$   acts  non-trivially on
the set of homotopy classes of strings.   Consider the permutation $\sigma=(134)(2)$ on the set $\{1,2,3,4\}$ sending 1 to
3, 3 to 4, 4 to 1, and 2 to 2.  Drawing the  string $ \alpha_\sigma$ we obtain that
$\overline {\alpha_\sigma}=\alpha_\tau$ where $\tau$ is the permutation $(124)(3)$.  The based matrices $T({\alpha_\sigma})$
and
$T(\alpha_\tau)$ can be explicitly computed. They are determined by the
following  skew-symmetric   matrices:  
$$
 \left [  \begin{array}{ccccc}  0& -2&0& -1& 3 \\
        2& 0& 1& 0& 3 \\
0& -1& 0& 0&2 \\
1& 0&0& 0&1 \\
-3& -3& -2& -1& 0 
\end{array} \right ],\,\,\,\,\,\,\,\,\,\,\,\,\,\,\,\,\,\,\,\, \left [  \begin{array}{ccccc}  0& -1&-2& 0& 3 \\
        1& 0& -1& 1& 3 \\
2& 1& 0& 1&2 \\
0&-1&-1& 0&1 \\
-3& -3& -2& -1& 0 
\end{array} \right ].$$
The based matrices   $T({\alpha_\sigma})$ and
$T(\alpha_\tau)$ are not isomorphic; this is clear for instance from the fact that the first matrix has a row with
three zeros while the second matrix does not have such a row. It is clear also that these based matrices are primitive. By Lemma
\ref{l:t51},  they are not homologous. Hence ${\alpha_\sigma}$ is not homotopic to
$\alpha_\tau=\overline
{\alpha_\sigma}$.

   \subsection{Remark}\label{fi:g5512}   For open strings one can define a refined version of based matrices incorporating
the splitting
$\arr= \arr^+ \cup \arr^- $.   A {\it refined based matrix}  is a  based matrix $(G, s, b)$ endowed with a spltting  of 
$G-\{s\}$ as a union of  disjoint subsets $G^+$ and
$G^-$. The moves  $M_1, M_2, M_3 $ on refined based matrices are defined  as above with $\overline G^{+}=G^+\cup \{g\},
\overline G^{-}=G^-$, 
$\tilde G^{+}=G^+, \tilde G^{-}=G^-\cup \{g\}$, and  $\hat G^{+}=G^+\cup \{g_1\},
\hat G^{-}=G^-\cup \{g_2\}$.  We leave further details to the reader.

   \section{Lie cobracket for strings}\label{cob}

   In this section we introduce a Lie cobracket in  the  free module generated by homotopy classes of strings.
This induces a Lie bracket in the module of homotopy invariants of strings.  

Throughout the section, we fix a 
commutative ring 
$R$ with unit.

		         \subsection{Lie   coalgebras}\label{su:61} 
We recall  here the  notion of a Lie coalgebra  dual to the one of a Lie algebra. To 
this end, we first  reformulate the notion of a Lie algebra.   For an $R$-module $L$,   denote by $\Perm_L$ the permutation 
$x\otimes 
y\mapsto y\otimes x$ in $L^{\otimes 2}=L\otimes  L$ and by $\tau_L$ the 
permutation $x\otimes y\otimes z\mapsto z\otimes x\otimes y$ in 
$L^{\otimes 
3}=L\otimes L \otimes   L$. Here and below  
$\otimes=\otimes_R$.  
 A Lie algebra over $R$ is an $R$-module $L$ endowed with an $R$-homomorphism (the 
Lie bracket)  $\theta:L^{\otimes 2}\to L$ such that   $\theta \circ 
\Perm_L=-\theta$ (antisymmetry) and    $$\theta \circ (\id_L\otimes 
\theta) 
\circ (\id_{L^{\otimes 3}}+\tau_L+\tau_L^2)=0 \in \Hom_R(L^{\otimes 3},L) $$
(the Jacobi identity).  Dually, a Lie coalgebra over $R$ is an $R$-module 
$A$ 
endowed with an $R$-homomorphism (the Lie cobracket)  $\nu:A\to 
A^{\otimes 2}  $ 
such that   $  \Perm_A \circ \nu=-\nu$   and  
			 \begin{equation}\label{jac}(\id_{A^{\otimes 
3}}+\tau_A+\tau_A^2)  \circ (\id_A\otimes \nu) \circ \nu=0 \in \Hom_R(A,  A^{\otimes 3}).
\end{equation}

A Lie 
coalgebra $(A,\nu)$ gives rise to the {\it dual Lie algebra}
  $A^*=\Hom_R(A,R)$ where the Lie bracket  $ A^*\otimes A^* \to 
A^*$ is  the homomorphism 
dual to $\nu$.  For $u,v\in A^*$, the value of $[u,v]\in A^*$ on     $x\in A$ is computed by
$$[u,v] (x)=\sum_i u(x^{(1)}_i) \, v (x^{(2)}_i)\in R$$
for any (finite) expansion $\nu(x)=\sum_i x^{(1)}_i \otimes x^{(2)}_i\in A\otimes A$.

A {\it homomorphism} of Lie coalgebras $(A,\nu)\to (A',\nu')$ is an $R$-linear 
homomorphism $\psi:A\to A'$ such that
$(\psi \otimes \psi) \nu (a)= \nu' \psi(a)$ for all $a\in A$. 
It is clear that the dual homomorphism $\psi^*:(A')^*\to A^*$ is a 
homomorphism 
of Lie algebras.

  \subsection{Lie   coalgebra of strings}\label{su71}
Let $\str$ be the 
set of 
homotopy classes of virtual strings and let $\str_0 \subset \str$ be its subset formed by the  homotopically non-trivial
classes.  Let
$\A_0=\A_0(R)$ be the free 
$R$-module freely generated by $\str_0$. We shall 
provide $\A_0$ with the 
structure 
  of a Lie coalgebra. 
  
 We begin with notation. For a  string $\alpha$,
let $\langle 
\alpha 
\rangle $ denote its class in $\str_0$ if $\alpha$ is homotopically non-trivial and set $\langle 
\alpha 
\rangle=0\in \A_0$ if $\alpha$ is homotopically trivial.   
 For an arrow $e=(a,b)$ of a string $\alpha$, denote
by $\alpha^1_e$ 
the string 
obtained from $\alpha$ by removing all arrows except
those with both 
endpoints 
in the interior of  the arc $ab$. (In particular, $e$ is removed.) Similarly, denote by $\alpha^2_e$ the string
obtained from 
$\alpha$ by  
removing all arrows except those with both endpoints
in the interior of  
$ba$. 
 Set
\begin{equation}\label{cobr}\nu (\langle \alpha
\rangle)= \sum_{e\in 
\arr(\alpha)}  \langle \alpha^1_e\rangle\otimes 
 \langle \alpha^2_e\rangle-\langle
\alpha^2_e\rangle\otimes 
 \langle \alpha^1_e\rangle  \in \A_0\otimes \A_0.\end{equation}

  \begin{lemma}\label{l:72} The $R$-linear homomorphism
$\A_0\to \A_0\otimes 
\A_0$ given 
on the generators of $\A_0$ by Formula \ref{cobr} is a
well-defined Lie 
cobracket.
		   \end{lemma}
                     \begin{proof} To show that $\nu$
is well-defined 
we must 
verify that   $\nu (\langle \alpha
\rangle)$ does not 
change under 
the homotopy moves (a$)_s$, (b$)_s$, (c$)_s$ on
$\alpha$. The  arrow 
added by 
(a$)_s$ contributes $0$ to the cobracket by the definition of    $\langle...\rangle$. The 
contribution of all the other arrows is preserved.
Similarly, the two 
arrows 
added by (b$)_s$ contribute  opposite terms to the cobracket  which is therefore preserved. 
Under (c$)_s$, all 
arrows 
contribute  the same before and after the move.	 
	
 The equality $\Perm_{\A_0} \circ \nu=-\nu$ is
obvious. We now 
verify Formula \ref{jac}. Let $\alpha$ be a string
with core circle 
$S$. We can expand   $(\id \otimes \nu)
( \nu (\langle 
\alpha 
\rangle))$ as a sum of   expressions $z(e,f)$
associated with     
ordered 
pairs of unlinked arrows $e,f\in \arr(\alpha)$.    
Note that the 
endpoints of 
$e,f$ split $S$ into four arcs meeting only at their
endpoints. The endpoints of $e$ (resp. $f$) bound one of these arcs, say $x$ (resp. $y$).   
 The 
other two arcs form $ S- (x\cup y)$ and lie \lq\lq
between" $e$ and 
$f$. Denote 
by $  \beta $ 	     (resp. $ \gamma $, $ \delta $) the
string obtained 
from 
$\alpha$ by removing all arrows except those with both
endpoints in the 
interior 
of $x$ (resp. of $y$, of   $ S- (x\cup y)$).  Set
$\varepsilon=+1$ if 
$e$ and 
$f$ are co-oriented, i.e., if their tails bound a
component of  $ S- 
(x\cup y)$. 
It is easy to see that
		     $$z(e,f)=  \varepsilon  (\langle
\beta\rangle \otimes 
\langle \delta\rangle \otimes  \langle \gamma\rangle 
- \langle 
\beta\rangle \otimes  \langle \gamma\rangle \otimes
 \langle 
\delta\rangle ).$$
A direct computation using this formula gives
$$(\id_{{\A_0}^{\otimes 3}}+\tau_{\A_0}+\tau_{\A_0}^2)
(z(e,f)+z(f,e))=0.$$		 
Thus $   \id_{{\A_0}^{\otimes 3}}+\tau_{\A_0}+\tau_{\A_0}^2
$ annihilates  
  $(\id \otimes \nu) ( \nu (\langle \alpha \rangle))$.
Hence $\nu$ is a 
Lie 
cobracket.\end{proof}

        \subsection{Lie coalgebra $\A$ and Lie algebra $\A^*$.}\label{su74}  Let
$\A=\A (R)$ be the free 
$R$-module freely generated by $\str$.  Since $\str=\str_0 \cup \{\Phi\}$  where $\Phi\in \str$ is the homotopy class of a trivial
string,  $\A=\A_0 \oplus R \Phi$.  The Lie cobracket $\nu$ in  $\A_0$ extends to $\A$ by $\nu (\Phi)=0$.

The   Lie cobrackets  in $\A_0$  and $\A$ induce 
  Lie brackets 
in  $\A_0^* 
=\Hom_R (\A_0 , R)$ and  $\A^* 
=\Hom_R (\A , R)$, respectively.  
Examples below show  that these  Lie cobrackets  and   Lie brackets
 are non-zero.  Clearly, $\A^*=  \A_0^*  \oplus R$ where the Lie bracket in $R$ is zero. 

The elements of
$\A^*
$ bijectively correspond to maps $\str \to R$, i.e., to 
$R$-valued 
homotopy invariants of strings. Thus, such invariants
form a Lie 
algebra.

	\subsection{Examples.}\label{su99}   (1) If   $\rank \alpha \leq 6$, then   $\nu (\langle
\alpha
\rangle)=0$. This follows from the fact that any string     of rank $\leq 2$ is homotopically trivial.  

(2) For any $p,q\geq 1$, we have $\nu (\langle\alpha_{p,q}\rangle)=0$.

(3)  Consider the string 
$\alpha=\alpha_\sigma$ of rank 7 where $\sigma$ is the permutation $ (123) (4) (576)$ of the set $\{1,2,\ldots, 7\}$. 
It follows from the definitions that
 $\nu(\langle \alpha\rangle)= \langle \alpha_{1,2} \rangle\otimes \langle \alpha_{2,1} \rangle- \langle \alpha_{2,1} \rangle
\otimes \langle \alpha_{1,2} \rangle$.
 As we know,  $\alpha_{1,2}$ and $ \alpha_{2,1}$   are homotopically non-trivial  strings      representing  distinct
generators of $\A$. Hence $\nu(\langle \alpha\rangle)\neq 0$.

(4) In generalization of the previous example pick any integers $p,q,p',q'\geq 1$ such that $p+ q\geq 3, p'+q'\geq 3$.
 Consider the string 
$\alpha=\alpha_\sigma$ of rank $m=p+q+p'+q'+1$  where $\sigma$ is the permutation   of the set
$\{1,2,\ldots, m\}$ defined by
$$\sigma(i)=\left\{\begin{array}{ll}
i+q ,~ {\rm {if}} 
\,\,\, 
1\leq i\leq p \\
\noalign{\smallskip}
i-p
,~ 
{\rm
{if}} \,\,\, p<i\leq p+q
\\
\noalign{\smallskip}
i,~ 
{\rm
{if}} \,\,\, i=p+q+1
\\
\noalign{\smallskip}
i+q'
,~ 
{\rm
{if}} \,\,\, p+q+1<i\leq p+q+1+p'  
\\
\noalign{\smallskip}
i-p'
,~ 
{\rm
{if}} \,\,\, p+q+1+p'  <i\leq m.
\end{array} \right.$$
  It follows from the definitions that
$$\nu(\langle \alpha\rangle)=  \langle \alpha_{p',q'} \rangle \otimes  \langle \alpha_{p,q} \rangle -  \langle \alpha_{p,q}
\rangle \otimes  \langle \alpha_{p',q'} \rangle .$$
Clearly,  $\nu(\langle \alpha\rangle)\neq 0$ unless $p= p'$ and  $q= q'$.

(5) Consider the numerical
invariants $u_1, u_2,\ldots \in  \A^* $ constructed in Section \ref{sn:g21}.  For  $p,p'\geq 1$, we  compute the
value of
$[u_p,u_{p'}]\in 
\A^*
$ on the  string $\alpha=\alpha(p,p',q,q')  $ defined in the previous example. Assume for  concreteness that   the numbers 
$p,p',q,q' $ are pairwise distinct. Then
$$[u_p,u_{p'}] (  \alpha )=u_p(\alpha_{p',q'} ) \,u_{p'} (\alpha_{p,q})- u_p(\alpha_{p,q} )\,  u_{p'} (\alpha_{p',q'}) =0-
(-q) (-q')=-qq'.$$
Hence $[u_p,u_{p'}]\neq 0$ for $p\neq p'$.

  \subsection{Filtration of $\A_0$.}\label{su735} Assigning to a string its homotopy rang and homotopy genus (see  
Section \ref{sn:g133}) we obtain two maps $hr, hg:\str_0 \to \ZZ$.  For   $r, g\geq 0$, set $$\str_{r,g}=\{\alpha \in
\str_0 \,
\vert \, hg(\alpha)\leq r,\,\,\,\, hg(\alpha)\leq g\}.$$
This set is finite since there is only a finite number of strings of rank $\leq r$.  The set  $\str_{r,g}$ generates a 
submodule  of $\A_0$ denoted $\A_{r,g}$.  This submodule is a free $R$-module of   rank $\#(\str_{r,g})$. 
Clearly, 
\begin{equation}\label{fil}\nu ( \A_{r,g})\subset \bigoplus_{p,q\geq 0, p+q<r} \A_{p,g} \otimes \A_{q,g} \subset  
\A_{r,g} \otimes \A_{r,g}.\end{equation}
Thus, each $\A_{r,g}$ a  Lie coalgebra.  The   inclusions $\A_{r,g}\hookrightarrow \A_{r',g'}$ for   $r\leq r', g\leq
g'$ make the family  $\{\A_{r,g}\}_{r,g}$ into a direct spectrum of Lie coalgebras. The equality $\A_0=\cup_{r,g}
\A_{r,g}$ shows that $\A_0=\injlim \{\A_{r,g}\}$.

The Lie cobracket in $\A_{r,g}$ induces
 a Lie bracket 
in   $\A_{r,g}^* 
=\Hom_R (\A_{r,g} , R)$.  Formula \ref{fil} implies that this Lie algebra is nilpotent. Restricting  maps  $\str_0 \to R$ to
$\str_{r,g}$ we obtain a Lie algebra homomorphism 
	$\A_0^*\to  \A_{r,g}^* $.  It is clear that $\A_0^*=\projlim \{\A^*_{r,g}\}$.

	\subsection{Relations with Lie coalgebras of
curves.}\label{su75}
	Let $\Sigma$ be a connected  surface and 
	$\hat \pi$ be the set of homotopy classes of closed
curves on $\Sigma$. 
(It can be identified with the set of conjugacy
classes in 
	$\pi=\pi_1(\Sigma)$.) 	 There is a  map $\psi : \hat \pi \to \str$ sending
each homotopy 
class 
of curves into the homotopy class of the underlying
strings.  Clearly, 
 $\psi  (\hat \pi )= \cup_r\,  \str_{r,g} $
where $g= g(\Sigma) $ is the genus of $\Sigma$.
Observe that the mapping class group of $\Sigma$ acts on $\hat \pi$ in the obvious way and  $\psi$ factors through the
projection of $\hat \pi$ to the set of orbits of this action.

Let $Z =Z(R)$ be the free
$R$-module with 
basis $\hat \pi$.  The map $\psi: \hat \pi \to \str$ induces an $R$-linear homomorphism $ Z\to \A$ whose image is equal to 
 $  \cup_r\,  \A_{r,g} $.  Composing this homomorphism with the   projection  $\A=\A_0\oplus R\Phi \to \A_0$  we obtain  an $R$-linear homomorphism  $\psi_0:Z\to \A_0$. 

The author defined in \cite{tu},
Section 8 a 
structure of a 
Lie coalgebra in $Z$. (In fact $Z$ is a Lie bialgebra,
but we shall not 
use the 
Lie bracket in $Z$.)     A direct comparison of the
definitions yields the following.

\begin{lemma}\label{l:75} The  map  $\psi_0 :Z\to \A_0$  
is a homomorphism 
of Lie 
coalgebras.  
		   \end{lemma}

Composing  $\psi_0$ with the inclusion $\A_0\hookrightarrow \A$ and dualizing  we obtain    a Lie algebra
 homomorphism 
 $\A^*\to Z^*$.

	\subsection{Applications.}\label{su101}  We claim that the product of strings is not commutative even up to homotopy: 
 there are  strings $\gamma,
\delta$ such that a product of $\gamma, \delta$ is not homotopic to a product of $\delta, \gamma$.   Consider the string
$\alpha$ constructed in Example
\ref{su99}.3.  Drawing a picture of $\alpha$, one  observes that $\alpha$ is a product of $\delta=\alpha_{2,1}$ with a
string,
$\gamma$, of rank 4 obtained  from $\alpha_{1,2}$ by  adding a \lq\lq small" arrow.  
Since $\gamma$ has a small arrow, it is easy to form a product of $\gamma$ with $\delta$ also having a
small arrow.  The resulting string, $\beta$,  is homotopic to a string of rank 6. Hence $\nu (\langle \beta\rangle )=0$
whereas  $\nu (\langle \alpha\rangle )\neq 0$. Therefore $\alpha$ is not homotopic to $\beta$.

\section{Virtual strings versus  virtual  knots}\label{vks}

 Virtual knots were introduced by L. Kauffman \cite{Ka} as a  
generalization of  classical knots.  We relate them to  virtual   strings by showing that each virtual knot gives rise to a
polynomial on virtual strings with coefficients in the   ring $\QQ [z]$.  As a technical tool, we
introduce a skein algebra of virtual knots and compute it in terms of  strings.

    \subsection{Virtual knots.}\label{vks:1}    We   define virtual knots  in terms of arrow diagrams following
\cite{gpv}.  An {\it arrow diagram} is a virtual string whose arrows are endowed
with   signs $\pm$.  By the core circle and the endpoints of  an arrow diagram, we mean the core circle and the endpoints of 
the   underlying   virtual string. 
The sign of an  arrow $e$ of an arrow diagram  is denoted
$\sign (e)$.   Homeomorphisms of arrow diagrams are defined as the homeomorphisms of the underlying strings preserving
the signs of all arrows.  The homeomorphism classes of       arrow diagrams will be  also called  arrow diagrams.  

We describe  
three moves (a$)_{ad}$, (b$)_{ad}$, (c$)_{ad}$ on arrow diagrams where $ad$ stands for
\lq\lq arrow diagram". Let $\alpha$ be an arrow diagram with core circle $S$. Pick two  distinct points 
$a,b \in S$  such that the (positively oriented) arc   $ab\subset S$ is disjoint from the set of   
endpoints of 
$\alpha$. The move (a$)_{ad}$   
adds to 
$\alpha$  the arrow $(a,b)$ with sign $+$ or $-$.   This move has two forms determined by the sign $\pm$.   The move
(b$)_{ad}$  acts on 
$\alpha$ as follows. Pick two   arcs   on $S$ disjoint from each 
other 
and   from the   endpoints of $\alpha$. Let $a,a'$ be the 
endpoints of 
the first arc (in an arbitrary order) and $b,b'$ be the endpoints of 
the second 
arc. The move adds 
to $\alpha$ two arrows $(a,b)$ and $(b',a')$ with opposite signs.  This move has eight forms depending on the
choice of the sign of $(a,b)$,  two possible choices for  $a$, and two possible choices for  $b$.   (This list of
eight forms of   (b$)_{ad}$ contains two equivalent pairs so that in fact the move  (b$)_{ad}$ has only six forms.)
 The
move (c$)_{ad}$     applies  to 
$\alpha$   when $\alpha$ has three arrows  with signs 
$((a^+,b),+), ( (b^+,c),+), ( (c^+,a),-)$ where  $a, a^+, b,b^+,c, c^+\in S$ such that
the   arcs   $aa^+$,   $bb^+$,     $cc^+$ are disjoint from each 
other and 
from the other   endpoints of $\alpha$. The   move (c$)_{ad}$  replaces 
these three  arrows    with the arrows  $ ((a,b^+),+), ( (b,c^+),+), ((c,a^+),-)$.  

By definition, a {\it virtual  knot}  is an equivalence class of arrow diagrams with respect to the equivalence relation
generated by  the  moves (a$)_{ad}$, (b$)_{ad}$,
(c$)_{ad}$ and  homeomorphisms. Note that our set of moves is somewhat different from the one in 
\cite{gpv} but  generates the same equivalence relation (cf. below). 

In the sequel the virtual knot represented by an arrow diagram $D$
will be denoted $[D]$.   A {\it trivial arrow diagram} having no arrows represents the {\it trivial virtual knot}. 

Forgetting the signs of arrows, we  can associate with  any arrow diagram $D$ its underlying virtual string $\underline
D$.   This induces   a \lq\lq forgetting" map  $K\mapsto \underline K$ from   the set of virtual knots 
   into the set of virtual strings. This map is surjective bur not  injective.  The theory of virtual knots is considerably reacher
than the theory of virtual strings. For instance,   the   fundamental  group of a virtual knot  \cite{Ka}  allows to
distinguish virtual knots with the same underlying strings.

   \subsection{From knots to virtual knots.}\label{vks:2}  Arrow diagrams are closely related  to the   standard
  knot diagrams on surfaces. An (oriented)
  knot diagram on an (oriented) surface $\Sigma$ is a (generic oriented)  closed curve on $\Sigma$ such that at
each its double point one of the branches of the curve passing through this point is distinguished. The distinguished branch is
called an {\it overcrossing} while the second branch passing through the same point is  called an {\it undercrossing}.  A knot diagram on
$\Sigma=\Sigma\times
\{0\}$ determines an (oriented)  knot in
$\Sigma\times
\RR$ by   pushing the overcrossings into  $\Sigma\times (0,\infty) $.  

Any knot diagram $d$ gives rise to an arrow diagram  $D(d)$  as follows.  First of all,  the closed curve underlying $d$ gives
rise to a virtual string, see Section \ref{sn:g12}. We provide each arrow  of this string with the sign of the corresponding
double point of $d$. This sign is $+$ (resp. $-$)  if the  pair  (a positive tangent vector to the
overcrossing 
 branch, a positive tangent vector to the undercrossing 
 branch) is positive  (resp. negative)  with respect to the orientation of $\Sigma$.  Our definition of the
arrow diagram associated with $d$ differs   from the one in \cite{gpv}: their arrow diagram  is obtained  from   ours by
reversing    all arrows  with  sign $-$.

There is a canonical mapping from
the set of isotopy classes of (oriented) knots in $\Sigma\times
\RR$ into the set of virtual knots.   It  assigns to a   knot    $K\subset \Sigma\times
\RR$ the virtual knot $[D(d)]$ where $d$  is a   knot diagram  on
$\Sigma$  presenting a knot   in  $\Sigma\times \RR$  isotopic to $K$.    The virtual knot  $[D(d)]$ does not depend  on 
  the choice  of $d$.  This follows from the fact that two  knot diagrams on
$\Sigma$  presenting isotopic knots  in  $\Sigma\times \RR$   can be obtained from each other  by ambient
isotopy in
$\Sigma$ and the     Reidemeister moves. Recall the standard list of the   Reidemeister
moves:  (1)  a  move adding a twist on the right (resp. left) of a branch; (2)  a move  pushing a
branch over another  branch and creating two crossings;  (3) a move pushing a branch over a crossing.  This list is
redundant. In particular,  the left move of type (1)    can be presented as a composition of  type (2) moves  and the inverse
to a right move of type (1).  One  move  of type (3)  together with moves of
type  (2)  is sufficient  to generate all moves of type (3) corresponding to various orientations on the branches  (see, for
instance, 
\cite{turr}, pp. 543--544).   As the generating move of type (3) we take the move  (c${)}^{-}$ described in Section
\ref{sn:g13}.   It remains   to observe that  the moves     (a$)_{ad}$, (b$)_{ad}$,
(c$)_{ad}$ on arrow diagrams are  exactly the moves induced by the right Reidemeister moves of type (1),  the  
 Reidemeister moves of type (2), and the move (c${)}^{-}$.
		     
 \subsection{Skein algebra of virtual knots.}\label{vks:3}    Let 
$R=\QQ[z]$ be the ring of polynomials in one variable $z$ with  rational coefficients.   Consider the polynomial  algebra $
R[\mathcal  K]$ generated by 
  the set of virtual knots $
\mathcal  K$.  This is a commutative associative algebra with  unit whose elements are polynomials in elements of $\mathcal  K$ with coefficients in
$R$.  
 We now introduce certain elements of $R[\mathcal  K]$ called \lq\lq skein relations".

   Pick an arrow diagram $D$ with core circle $S$ and pick an arrow $e=(a,b)$ of $D$ with sign $+$ (here $a,b\in S$).  Let
$D^-_e$ be the same arrow diagram with the sign of $e$ changed to $-$.  Let $D'_e$ be the arrow diagram obtained from
$D$ by removing all arrows with at least one endpoint on the arc $ba\subset S$. Let $D''_e$ be the arrow diagram obtained
from
$D$ by removing all arrows with at least one endpoint on the arc $ab\subset S$.  The skein relation corresponding to
$(D,e)$ is
 $[D]-[D^-_e]-z [D'_e] [D''_e]\in R[\mathcal  K]$.

The ideal of  the algebra $ R[\mathcal  K]$ generated by the trivial virtual  knot  and the skein relations
(determined by all the  pairs
$(D,e)$ as above) is called the {\it skein ideal}.   The quotient of  
$ R[\mathcal  K]$  by this ideal   is called the {\it skein algebra of virtual knots} and denoted  $\mathcal {E}$. The  next
theorem computes
$\mathcal {E}$ in terms of strings. Recall the set 
 $\str_0$   of  non-trivial 
homotopy classes of virtual strings, cf. Section \ref{su71}.

\begin{theor}\label{th:eee}
                    There is a canonical $R$-algebra isomorphism $\nabla: \mathcal {E} \to R[\str_0]$ where $R[\str_0]$ is the
polynomial algebra  generated by 
 $
\str_0$.
                     \end{theor}

This  theorem allows us to associate with any virtual knot $K$ a polynomial $\nabla(K)\in R[\str_0]$.   It will be clear from
the definitions that
$$\nabla(K)= \langle \underline K \rangle+\sum_{n\geq 2} z^{n-1}  \,\nabla_n (K)$$
where $ \nabla_n (K) $ is a homogeneous element of  $\QQ [\str_0]$ of degree $n$ which is non-zero only for a finite set of
$n$.  Combining
$\nabla$ with homotopy invariants of strings we obtain  invariants of virtual knots.  For example, composing $\nabla$ with the
algebra homomorphism
$R[\str_0] \to R[t]$ sending the homotopy class of a string $\alpha$ into the polynomial $u(\alpha)(t) $, we obtain an
algebra homomorphism  $\mathcal {E} \to R[t]=\QQ[z,t]$. This gives a   2-variable polynomial invariant of virtual knots. 

The constructions above can be applied to virtual knots derived from geometric knots in  (a surface)$\times \RR$ as in
Section
\ref{vks:2}. The resulting invariants are interesting only in the case when the genus of the surface is
at least $2$.  This is due to the fact that the strings realized by curves on a surface  of genus 0 or 1 are homotopically
trivial.

Theorem \ref{th:eee} will be proven in the next section.  Here we   give   an explicit expression for the value of  $\nabla$
on the generator  
$[D]\in \mathcal {E}$ represented by an arrow diagram
$D$. We need a few definitions. 
  The endpoints of  the
arrows of  $D$ split the core circle of $D$ into  (oriented) arcs called the {\it edges} of $D$.  Denote the set of  
edges of $D$ by
$\edg(D)$.   Each endpoint  $a$ of an arrow of $D$ is adjacent to two edges  $a_-,a_+\in \edg(D)$, respectively
incoming and   outgoing with respect to $a$. 
For an integer $n\geq 1$, an {\it $n$-labeling} of $D$ is a   map  $f:\edg(D)\to \{1,2,...,n\}$ satisfying the
following condition: for any  arrow $e=(a,b)$ of $D$,  either 

 (i) $f(a_+)=f(a_-), 
f(b_+)=f(b_-)$ or 

 (ii) 
$f(a_+)= f(b_-) \neq  f(a_-)= f(b_+)$ and  $\sign( f(a_-)-f(a_+))=\sign (e)$.  

The arrows $e$ as in (ii) are called {\it $f$-cutting arrows}.  The number of   $f$-cutting
arrows 
  of 
$D$ is denoted $\vert f\vert$ and the   number of $f$-cutting
arrows 
  of 
$D$ with $\sign =-1$  is denoted $\vert f\vert_-$.
Note that the value of $f$ on two adjacent 
edges
$a_-,a_+\in
\edg(D)$ may differ only when 
$a$ is  an endpoint of an  $f$-cutting arrow.  Therefore   $\vert f\vert \geq \# f(\edg
(D))-1$.  For $i=1,...,n$,  let ${\underline D}_{f,i}$   be the string  obtained from $D$ by removing all arrows except  the
arrows
$ (a,b)$ with $f(a_+)=f(a_-)= 
f(b_+)=f(b_-)=i$ (and forgetting the signs of the arrows).    

Let $\lbl_n (D)$ be the set of $n$-labelings $f$
of
$D$ such that   $f(\edg
(D))=\{1,...,n\}$,   $\vert f\vert=n-1$, and the  $f$-cutting arrows of $D$ are  pairwise unlinked
(in the sense of Section 
\ref{sn:g21}).  Then
\begin{equation}\label{etae}\nabla([D])=\sum_{n=1} \sum_{f\in  \lbl_n (D)} \frac {(-1)^{\vert f\vert_-} z^{n-1}}{n!}
\prod_{i=1}^n
\langle {\underline D}_{f,i}  
\rangle  \in R[\str_0].\end{equation}
The expression on the right-hand side is finite since $\lbl_n (D)=\emptyset$ for $n>\# \edg(D)$.  The set
 $\lbl_1 (D)$ consists of only one element $f=1$ so that  the   free term of
$\nabla([D])$  is   $\langle {\underline
D} 
\rangle$.

\section{Proof of Theorem \ref{th:eee}}\label{vksss}

The proof of Theorem \ref{th:eee}  largely follows the proof of Theorems 9.2 and 13.2 in \cite{tu}.  We therefore expose
only the main lines of the proof.  The key point  behind Theorem \ref{th:eee} is the existence of a natural comultiplication in
$\mathcal E$ and we define it first.  Then we   construct $\nabla$ and   prove that it is an isomorphism.

 \subsection{Comultiplication in $\mathcal E$.}\label{vks:4} We need   to
study more extensively the  labelings of   arrow diagrams defined at the end of the previous section.
Let $D$ be an arrow diagram with core circle $S$.  Each  $n$-labeling $f$ of   
$D$  gives rise to  $n$ 
monomials  $D_{f,1}, ..., D_{f,n} \in \mathcal
{E}$  as follows.    Identifying  $a=b$ for every $f$-cutting arrow  $(a,b)$ of $D$, we
transform  
$S$ into a 4-valent graph,  $\Gamma^f$, with $\vert f\vert$  vertices.   The  projection $S\to
\Gamma^f$  maps the non-$f$-cutting arrows of
$D$ into \lq\lq arrows" on  $\Gamma^f$, i.e.,  into ordered pairs  of (distinct)  generic points of $\Gamma^f$. The
labeling
$f$   induces a labeling of the edges of
$\Gamma^f$ by the numbers $1,2,...,n$.  It follows from the definition of  a labeling that for each $i=1,...,n$, the
union of edges of
$\Gamma^f$ labeled with
$i$ is  a disjoint union of $r_i=r_i(f)\geq 0$   circles   $S^i_1,..., S^i_{r_i}$. 
The orientation of   $S$ induces an
orientation of the edges of $\Gamma_f$ and of these circles.  We  transform each  
circle  $S^i_j$ with $j=1,..., r_i$  into an arrow diagram by adding to it all the arrows of $\Gamma^f$ with both endpoints
on
$S^i_j$.  The signs of these arrows are by definition the signs of the corresponding non-$f$-cutting arrows of $D$.   
Set $$D_{f,i}= \prod_{j=1}^{r_i} \, [S^i_j] \in \mathcal E.$$

For    any $n\geq  2$, denote  
$\Lbl_n  (D)$ the set of
$n$-labelings
$f$ of
$D$ such that  the $f$-cutting arrows of $D$ are  pairwise unlinked.  The latter condition can be reformulated in terms of
the numbers $r_1(f),...,r_n(f)$ introduced above: $f\in \Lbl_n  (D)$ if and only if $r_1(f)+ \cdots + r_n (f)= \vert f\vert
+1$.   For
$f\in \Lbl_n (D)$, set  
$$\Delta(D,f)= (-1)^{\vert f\vert_-} z^{\vert f\vert}\,  D_{f,1}  \otimes   D_{f,2}  \otimes \cdots  \otimes D_{f,n} \in
\mathcal {E}^{\otimes n}$$
where $\mathcal {E}^{\otimes
n}$ is the tensor product over $R$ of $n$ copies of $\mathcal {E}$.

  By a comultiplication in $\mathcal E$, we mean a coassociative algebra homomorphism $\Delta:\mathcal
{E}\to
\mathcal {E}{\otimes }\mathcal {E}$. (The coassociativity means that $(\id \otimes \Delta) \Delta= (\Delta \otimes \id
)\Delta$.)   We claim that the formula 
$$\Delta ([D])=   \sum_{f\in  \Lbl_2 (D)} \Delta(D,f)\in \mathcal {E}{\otimes
}\mathcal {E}$$
   extends  by multiplicativity to a well-defined  comultiplication    in $\mathcal E$. 
This can be deduced from     \cite{tu},  Theorem  9.2  or  proven directly repeating the same arguments.  We explain
how to deduce our claim from    \cite{tu}.  Comparing the  definition of 
$\Delta ([D])$ with the comultiplication in the algebra of skein classes of knots in (a surface)$\times \RR$ given in  \cite{tu},
we observe that they correspond to each other provided $D$  underlies a knot diagram on the surface.  (The variables
$h=h_1, \overline h=h_{-1}$ used in  \cite{tu} should be replaced with $0$ and $z$, respectively.  After the substitution 
$h=0$,  we can consider only labelings  satifying - in the notation of  \cite{tu} - the condition $\vert \vert f \vert \vert =-
\vert f\vert$ which  translates here as the assumption that   the  $f$-cutting arrows of $D$ are  pairwise unlinked.)   
The  results of 
 \cite{tu} imply that  if
  a  move 
(a$)_{ad}$, (b$)_{ad}$, (c$)_{ad}$  on
$D$ underlies a Reidemeister move on a knot diagram, then  $\Delta ([D])$ is preserved under  this move. Since any
arrow diagram $D$ underlies a knot diagram on a surface  and any move  (a$)_{ad}$, (b$)_{ad}$, (c$)_{ad}$  on
$D$ can be induced  by a  Reidemeister move, we conclude that  $\Delta ([D])$ is 
invariant under the moves  (a$)_{ad}$, (b$)_{ad}$, (c$)_{ad}$  on
$D$.  Therefore the formula $[D] \mapsto \Delta([D])$  yields a well-defined mapping $\mathcal K \to     \mathcal
{E}{\otimes }\mathcal {E}$. This mapping  uniquely extends    to an  algebra
homomorphism 
$ R[\mathcal
{K}] \to
\mathcal {E}{\otimes }\mathcal {E}$. 
The  results
of 
 \cite{tu} imply that  for  an arrow diagram   $D$  underlying a knot diagram  on a surface and  any arrow $e$ of $D$ with
$\sign (e)=+$, the   skein relation  
 $  [D]- [D^-_e]-z  [D'_e]  \,  [D''_e]  $ lies in the kernel of the latter homomorphism.   The condition that  
$D$  underlies a knot diagram    is verified for all $D$.   Therefore the  homomorphism $ R[\mathcal
{K}] \to
\mathcal {E}{\otimes }\mathcal {E}$ annihilates the skein ideal and induces  an  algebra
homomorphism 
$ \Delta: \mathcal
{E}\to
\mathcal {E}{\otimes }\mathcal {E}$.  
The coassociativity of $\Delta$  follows from the easy 
formulas
$$(\id \otimes \Delta) \Delta ([D])=  \sum_{f\in  \Lbl_3 (D)} \Delta(D,f)= (\Delta \otimes \id
)\Delta ([D])$$
 (cf.   \cite{tu}, p. 665).  More generally, 
   for any $n\geq 2$, the value  on $[D] \in \mathcal E$ of the iterated homomorphism 
$$\Delta^{(n)}= (\id^{\otimes (n-1)} \otimes \Delta)\circ  (\id^{\otimes (n-2)} \otimes \Delta) \circ ...\circ 
(\id  \otimes \Delta) \Delta: \mathcal {E}\to \mathcal {E}^{\otimes
 (n+1)}$$
  is computed by
$$\Delta^{(n)} ([D])= \sum_{f\in  \Lbl_{n+1} (D)} \Delta(D,f).$$
Note for the record
that   each arrow diagram $D$ admits   constant 2-labelings $f_1, f_2$   taking values $ 1,2$ on all edges,  respectively. 
The corresponding  summands of $\Delta ([D])$ are 
$\Delta(D,f_1)=[D] \otimes 1$ and $\Delta(D,f_2)= 1  \otimes [D]$.

  \subsection{Homomorphism $\nabla: \mathcal {E} \to R[\str_0]$.}\label{vks:6} There are two obvious $R$-linear
homomorphisms  $\varepsilon:  \mathcal {E}\to  R$ and  $ q:  \mathcal {E} \to R[\str_0]$.   The homomorphism
$\varepsilon$ sends
$1\in  \mathcal {E}$ into $1\in R$ and sends all   virtual knots and their  non-void products into $0$. 
 The homomorphism
$q$ sends
  $1$ and  all  products of $\geq 2$ virtual knots into
$0$  and sends  a virtual knot $K$ into $\langle \underline K \rangle$.    Tensorizing $q$ with itself,  we obtain for all 
$n\geq 1$ a homomorphism $q^{\otimes n}: \mathcal {E}^{\otimes n} \to R[\str_0]^{\otimes n}$.  Let
$s_n: R[\str_0]^{\otimes n}
\to  R[\str_0]$ be the $R$-linear homomorphism sending $a_1\otimes \cdots \otimes a_n$ into    $(n!)^{-1} a_1\cdots
a_n$. Set 
$$\nabla=\varepsilon + q+ \sum_{n\geq 2}  \, s_n \,q^{\otimes n} \Delta^{(n-1)}:  \mathcal {E} \to R[\str_0] $$
where  $\Delta^{(1)}=\Delta$.  It is clear that $\nabla$ is  $R$-linear. The same  argument as in   \cite{tu},  Lemma 13.4
shows that $\nabla$ is an algebra homomorphism.  Computing $\nabla$ on
the skein class of an arrow diagram $D$, we obtain 
$$\nabla([D])=\sum_{n\geq 1}   \, s_n \,q^{\otimes n}   \sum_{f\in  \Lbl_{n} (D)} \Delta(D,f)= \sum_{n\geq 1}       \sum_{f\in  \Lbl_{n} (D)} 
\frac {(-1)^{\vert f\vert_-} z^{\vert f\vert}}{n!} \,  \prod_{i=1}^n  q(D_{f,i}).$$
Note that  $q(D_{f,i})=0$ unless $r_i(f)=1$ in which case  $q(D_{f,i})=\langle \underline D_{f,i}\rangle$. 
For a labeling $f\in \Lbl_n(D)$ the equalities $r_1(f)=...=r_n(f)=1$ are equivalent to the inclusion $f\in \lbl_n(D)$.  
This  yields  Formula \ref{etae}.

Observe that   $\nabla([D])$ is a sum of   $\langle \underline
D\rangle$ and a polynomial  in strings  of rank $<\rank D$.  An  induction on the rank of strings shows that the image of
$\nabla$ contains all strings. Therefore $\nabla$ is surjective. 

 The proof of the injectivity of $\nabla$  is based on the
following   lemma. 

\begin{lemma}\label{l:jae} There is a $\QQ$-valued function  $\eta$ on the set of isomorphism classes of (finite) oriented
trees such that the following holds:

(i) if $T$ is  a tree with one vertex and no edges, then $\eta(T)=1$;

(ii) if an oriented tree $T'$ (resp. $U$) is obtained from an oriented tree $T$ by reversing the orientation of an edge
$e$ (resp. by contracting $e$ into a point), then $\eta(T)+\eta(T')+\eta(U)=0$;

(iii)   if an oriented tree $T'$ (resp. $T''$) is obtained from an oriented tree $T$ by replacing two distinct edges with
common origin $ab, ac$ by $ab, bc$ (resp. by $ac, cb$) and if $U$ is obtained from $T$ by identifying $b$ with $c$ and
$ab$ with $ac$, then $\eta(T)=\eta(T')+ \eta(T'')+\eta(U)$.
\end{lemma}

In this lemma  by an edge $ab$ we mean an {\it oriented} edge directed from $a$ to $b$.

Lemma \ref{l:jae} was first established in \cite{tu}, Theorem 14.1  where it is also shown that   $\eta$  is unique 
(we shall not need this).  The construction   in \cite{tu} is indirect and does not provide an explicit  formula for $\eta$.  Such
a formula was pointed out by Fran\c cois Jaeger  \cite{ja}.  The following proof of Lemma \ref{l:jae}  is a simplified version
of the proof given by   Jaeger 
\cite{ja}.

                     \begin{proof} By a {\it forest} we shall mean  a disjoint union of a finite family of  finite oriented trees.  The set
of vertices of a forest $F$ is denoted $V(F)$. For a forest
$F$ and an integer $n\geq 1$,  denote by  $C_n(F)$ the set of surjective mappings $f:V(F) \to \{1,...,n\}$ such that for
every edge $ab$ of $F$ we have $f(a)<f(b)$.  This set is empty for  $n> \#(V(F))$.  Set
$$\eta(F)= \sum_{n\geq 1} \frac {(-1)^{n+1}} {n}\, \# (C_n(F)) \in \QQ.$$
We claim that $\eta$ satisfies all the conditions of the lemma.  Condition (i) is obvious. Condition (iii) is a  direct  corollary of
the definitions. Indeed for all $n$,  the set $C_n(T)$ splits as a disjoint union of the sets $C_n(T'), C_n(T''), C_n(U)$. Hence
$\# (C_n(T))=\# (C_n(T'))+\# (C_n(T''))+\# (C_n(U))$ and   $\eta(T)=\eta(T')+ \eta(T'')+\eta(U)$.  It remains to verify 
(ii).   Let $F$ be obtained from $T$ by removing the interior of the edge $e$.  For all $n$,   
the set $C_n(F)$ splits as a disjoint union of the sets $C_n(T), C_n(T'), C_n(U)$. Hence 
$\# (C_n(F))=\# (C_n(T))+\# (C_n(T'))+\# (C_n(U))$ and  $\eta(F)=\eta(T)+ \eta(T')+\eta(U)$. Thus we need only to
prove that  $\eta(F)=0$ for every forest $F$ with two components $T_1,T_2$.

For  non-negative integers $n,k_1,k_2 $, denote by $C_n(k_1,k_2)$ the set of pairs $(l_1,l_2)$ where for
$i=1,2$, $l_i$ is an order-preserving
injection from $\{1,..., k_i\}$ into $\{1,...,n\}$ and   
$l_1(\{1,..., k_1\})\cup l_2(\{1,..., k_2\})= \{1,..., n\}$. Having $g_1\in C_{k_1} (T_1)$, $g_2\in C_{k_2} (T_2)$
and having  $(l_1,l_2)\in C_n (k_1,k_2)$ we   define a mapping  $f=f(g_1,g_2,l_1,l_2):V(F) \to \{1,...,n\}$ by
  $f(v)=l_1 g_1(v)$ for $v\in V(T_1)$ and  $f(v)=l_2 g_2(v)$ for $v\in V(T_2)$.
Clearly,  $f\in C_n(F)$.  It is obvious that any      $f\in C_n(F)$ can be uniquely presented in the form 
$f=f(g_1,g_2,l_1,l_2)$ where $g_i\in C_{k_i} (T_i)$ with $k_i=\#( f(V(T_i)))\geq 1 $ for $i=1,2$.  Therefore
$$\eta(F)=  \sum_{n\geq 1}   \frac
{(-1)^{n+1}} {n}\, \left ( \sum_{k_1,k_2 \geq 1} \,\sum_{g_1\in C_{k_1} (T_1),g_2\in C_{k_2} (T_2)}  \# (C_n
(k_1,k_2)) \right )$$
$$ =
\sum_{k_1,k_2 \geq 1} \# (C_{k_1} (T_1)) \, \# (C_{k_2} (T_2)) 
\sum_{n\geq 1}   \frac {(-1)^{n+1}} {n}\, \# (C_n (k_1,k_2)) .$$
Thus  it is enough  to prove that for all $k_1 \geq 1,k_2\geq 1$, the numbers $c_n (k_1,k_2)=\# (C_n (k_1,k_2))$ verify
\begin{equation}\label{fra} \sum_{n\geq 1}   \frac {(-1)^{n+1}} {n}\, c_n (k_1,k_2)=0.\end{equation}
Clearly, $c_n (k_1,k_2)$ is the number of pairs $(S_1,S_2)$ where $S_1,S_2$ are subsets of $\{1,...,n\}$ such that 
$S_1\cup S_2=\{1,...,n\}$, $\#(S_1)=k_1$, $\#(S_2)=k_2$. In particular, $c_n (k_1,k_2)=0$  if $k_1+k_2<n$ or  
$k_1>n$ or $k_2>n$.   For any
$n\geq 1$ and commuting variables
$x, y $, 
$$ (x+y+xy)^n=\sum_{  k_1,k_2\geq 0 } c_n (k_1,k_2)\, x^{k_1} y^{k_2}.$$
Therefore  
$$\log (1+x+y+xy)= \sum_{n\geq 1} \frac{ (-1)^{n+1} }{n} (x+y+xy)^n
$$
$$= \sum_{n\geq 1} \,\,\sum_{  k_1,k_2\geq 0 }  \frac {(-1)^{n+1}}{n} c_n (k_1,k_2) \, x^{k_1} y^{k_2}
= \sum_{  k_1,k_2\geq 0 } \,\,    \left (\sum_{n\geq 1} \frac {(-1)^{n+1}}{n} c_n (k_1,k_2) \right )\, x^{k_1} y^{k_2}.$$
Since
$$ \log (1+x+y+xy) =\log ((1+x)(1+y))=\log (1+x) +\log (1+y),$$
the terms with $k_1 \geq 1,k_2\geq 1$ in the above series must vanish.  This gives  Formula \ref{fra}. \end{proof} 

\subsection{The injectivity of $\nabla: \mathcal {E} \to R[\str_0]$.}\label{vks:624} We begin by associating with any
virtual string $\alpha$ an element $\zeta (\alpha)\in  \mathcal {E}$.  Let $S$ be the core circle of $\alpha$. 
A {\it  surgery} along an   arrow  $(a,b)\in \arr(\alpha)$ consists in picking two    (positively oriented) arcs    
$aa^+, bb^+\subset S$ and  then quotienting the  complement  of their interiors  $S-((aa^+)^\circ \cup (bb^+)^\circ)$ by
$a=b^+, b=a^+$.  It is understood that the arcs    
$aa^+, bb^+ $ are small enough not to    contain  endpoints of $\alpha$ besides $a,b$, respectively.  Such a surgery
transforms
$S$ into  two disjoint oriented circles.  We make each of them into a string by adding  all the arrows of
$\alpha$ with both endpoints  on the arc
$a^+b$ (resp. on $ba^+$).  (The arrows of $\alpha$ with one endpoint on $ab$ and the other one on $ba$
disappear under  surgery.)

Let us call a
set
$F\subset
\arr(\alpha)$ {\it special} if the arrows of $\alpha$ belonging to $F$ are pairwise unlinked.  
Applying surgery inductively to all arrows of
$\alpha$ belonging to a special set 
$F$, we transform
$\alpha$ into    $n=\#(F)+1$ strings.  Providing all the arrows of these   strings
with sign $+$, we obtain $n$ arrow diagrams    $D^F_1,..., D^F_n$.  Note that  
 they
  have together  at most     $\#(\arr(\alpha))-\#(F) $ arrows. We now  define an oriented graph
$\Gamma_F$. The vertices of $\Gamma_F$ are the symbols   $v_1,..., v_n$.  Two verices $v_i, v_j$
are related by an oriented edge   leading from $v_i$ to $v_j$ if there is an   arrow  $(a,b)\in F$ such that
the arcs $aa^+, bb^+\subset S$ involved in the surgery along this arrow lie  on the core circles of  $D^F_i, D^F_j$,
respectively. It is easy to see that $\Gamma_F$ is a tree. Set
$$\zeta(\alpha)=\sum_{F\subset \arr(\alpha)} \eta(\Gamma_F)\, z^{\#(F)} \prod_{i=1}^{\#(F)+1} [ D^F_i]\in \mathcal
{E}$$ where $F$ runs over all special subsets of $\arr(\alpha)$.  The summand corresponding to $F=\emptyset$ is the
string
$\alpha$ itself with sign  $+$ on all arrows.  
 
The key property of $\zeta(\alpha)  \in \mathcal
{E}$ is   its invariance under the   basic homotopy moves on $\alpha$. This follows from  \cite{tu}, Lemma 15.1.1
in the case where the moves are realized geometrically by   homotopy of a curve   realizing $\alpha$  on a surface.   Since
the homotopy moves can be always realized geometrically, the result follows. 
The mapping $\alpha\mapsto \zeta(\alpha)$  extends by multiplicativity to an algebra homomorphism $
R[\str_0]\to \mathcal
{E}$ denoted also  $\zeta$.

We can now prove the injectivity of $\nabla$. For $r\geq 0$, denote by $B_r$ the $R$-submodule of $\mathcal E$ 
additively generated  by  monomials $[D_1] [D_2] \cdots [D_n]$ such that the total number of arrows in the arrow
diagrams $D_1, D_2, ...,D_n$ is less than or equal to $r$.  Clearly, $0=B_0\subset B_1\subset ...$ and $\cup_r B_r=
\mathcal E$.  Pick $b= [D_1]
[D_2]
\cdots [D_n]\in B_r$.  Using the skein relation in $\mathcal E$ it is easy to see that    $b (\modu B_{r-1})\in
B_r/B_{r-1}$ does not depend on the signs of the arrows of $D_1,...,D_n$.  This observation, Formula \ref{etae}  and the
definition   of $\zeta $ imply  that 
$(\zeta
\nabla )(b) -b \in B_{r-1}$. Therefore $(\zeta \nabla-\id)^r(b) =0$. The inclusion $b\in \Ker \nabla$ would imply  
$b=0$. Thus  $B_r \cap \Ker \nabla  =0$.   Since $\cup_r B_r=
\mathcal E$, we obtain   $\Ker \nabla=0$. 

\subsection{Remarks.}\label{vks:624ddd}   The comultiplication   $\Delta$ defined  in Section \ref{vks:4} makes $\mathcal
E$ into a Hopf algebra over $R$. Its counit is the    homomorphism  $\varepsilon:  \mathcal {E}\to R$ 
 used in the definition of  $\nabla$.   For an arrow diagram $D$, denote by
$\tilde D$ the same diagram with opposite signs  on all
arrows.  The transformation $[D] \mapsto -[\tilde D]$ preserves the skein relation  and therefore induces an algebra
automorphism of
$\mathcal E$.  This automorphism is an antipode for $\mathcal
E$. This follows from the corresponding theorem for the skein algebras of curves on surfaces conjectured in \cite{tu} and
proven  in \cite{crr}  and independently in 
\cite{pr}.     In the construction of the Hopf algebra
$\mathcal E$ instead of the ground ring
$R=\QQ[z]$ we can use  $\ZZ[z]$.  It is only to construct the homomorphisms $\nabla$ and  $\zeta$   that
we need    $\QQ$.

 \section{Algebras and groups associated with strings}\label{fur}

  We discuss   here  various  algebraic structures associated with 
closed and open strings.    We begin by recalling the notion of
a spiral Lie coalgebra  and 
several related definitions from \cite{tu}, Section 11.    Throughout the section, $R$ is a 
commutative ring 
  with unit and $\otimes=\otimes_R$.

    \subsection{Spiral Lie coalgebras.}\label{su:6111}  For a Lie coalgebra $(A,\nu:A\to A^{\otimes 2})  $ over     $R$ and
an integer 
$n\geq 
1$, set  
   $$\nu^{(n)}=(\id_A^{\otimes (n-1)} \otimes \nu) \circ \cdots \circ 
  (\id_A  \otimes \nu) \circ \nu:A \to A^{\otimes (n+1)}.$$
In particular, $\nu^{(1)}=\nu$.
A  Lie coalgebra $(A,\nu)$ over $R$  is {\it spiral}, 
if $A$ is 
free as the $R$-module and the filtration $\Ker \nu^{(1)}\subset \Ker 
\nu^{(2)}\subset 
\cdots $ exhausts $A$, i.e., $A=\cup_{n\geq 1} \Ker \nu^{(n)}$. 

 Assume from now on that $A$ is spiral. 
The  dual Lie algebra $A^*=\Hom_R(A,R)$   has the following completeness property. 
Consider 
the lower central series $A^*=A^{*(1)}\supset A^{*(2)}\supset \cdots$ of 
$A^*$ 
where $A^{*(n+1)}= [A^{*(n)}, A^*]$ for $n\geq 1$. Let $a_1,a_2,\ldots 
\in A^*$ 
be an infinite sequence such that for any $n\geq 1$ all terms of the 
sequence 
starting from a certain place belong to $A^{*(n)}$. Clearly, if $x \in 
\Ker 
\nu^{(n)}$ and $a\in A^{*(n+1)}$, then $a(x)=0$.  Since $A=\cup_{n } 
\Ker 
\nu^{(n)}$, 
the sum $a(x)=a_1(x)+a_2(x)+\cdots$ contains only a finite number of 
non-zero 
terms for every $x\in A$.  Therefore $a(x)$ is a well-defined element of $R$. The formula $x\mapsto a(x):A\to R$   defines  
an element  of
$   A^* $  denoted $a_1+ a_2+ \cdots$ and called the (infinite) 
sum of 
$a_1,a_2,\ldots$.  A similar argument shows that 
 $\cap_n  A^{*(n)}=0$ and 
the natural Lie algebra homomorphism   $A^*\to \projlim_n (A^*/A^{*(n)})$ 
is an 
isomorphism.

  For   $a,b\in A^*$,    
consider the  sum
$$\mu (a,b)= a+b + \frac{1}{2} [a,b]+ \frac{1}{12} ([a,[a,b]] +[b, 
[b,a]])+ 
\cdots \in A^*$$
where the right-hand side is the Campbell-Hausdorff series for $\log 
(e^a e^b)$, 
see \cite{se}. The resulting mapping $\mu:A^*\times A^*\to A^*$ is a 
group 
multiplication in $A^*$. Here $a^{-1}=-a$ and $0$ is the group unit. 
The group  
$(A^*,\mu) $ is denoted $\Exp A^*$. Heuristically, this is 
the \lq\lq Lie group" with Lie algebra $A^*$. The equality  $A^*= \projlim_n (A^*/A^{*(n)})$ 
implies that   the  group $\Exp A^*$ 
is 
pro-nilpotent.

Consider the symmetric (commutative and 
associative) algebra of $A$:
$$S=S(A)=\oplus_{n\geq 0} S^n (A).$$   Here $S^0(A)=R$, $S^1(A)=A$, and $S^n(A)$  
is the 
$n$-th symmetric tensor power of $A$ for $n\geq 2$. The unit $1\in R= 
S^0(A)$ is 
the unit of $S$. The group multiplication $\mu:A^*\times A^*\to A^*$ 
  induces a comultiplication   $S\to S\otimes S$  as follows. Since 
$A$ is a 
free $R$-module, the natural map $A\to (A^*)^*$ extends to an embedding 
of $S$ 
into the algebra of $R$-valued functions on $A^*$. We can identify
 $S$ with the image of this embedding. Similarly, we can identify 
$S\otimes S$ 
with  an algebra of $R$-valued functions on $A^*\times A^*$.  It is 
easy to 
observe that for any $x\in S$, we have $x\circ \mu\in S\otimes S$. 
Indeed, it 
suffices to prove this for   $x\in A$. Then  
$x\in \Ker 
\nu^{(n)}$ for some $n$ so that $x$ annihilates all but finite number of terms of the  
Campbell-Hausdorff series. Our claim follows then from the duality 
between the 
Lie bracket in $A^*$ and the Lie cobracket $\nu$. For example, if $n=3$ 
and 
$\nu^{(2)}(x)=\sum_i \alpha_i\otimes \beta_i \otimes \gamma_i \in 
A^{\otimes 3}$, 
then
 $$x\circ \mu=x\otimes 1+1\otimes x+ \frac{1}{2} \nu(x)+\frac{1}{12}
 \sum_i (\alpha_i  \beta_i \otimes \gamma_i+ \gamma_i\otimes \alpha_i  
\beta_i )
.$$
 The formula $\Delta (x)=x\circ \mu$ defines a coassociative 
comultiplication in 
$S$. It has a  counit $S\to R$ defined as  the projection to $S^0(A)=R$. The antipode $S\to S$ is the 
algebra 
homomorphism sending any $x\in A$ into $-x\in A$.  A routine check 
shows that 
  $S$   is a  (commutative)
Hopf 
algebra.
 Heuristically, it should be viewed as the Hopf algebra of  $R$-valued functions on
the 
group $\Exp A^*$ or as the Hopf   dual of the universal enveloping   
algebra of 
  $A^*$. 
 
The construction of     $\Exp A^*$ and    $S(A)$ can be  generalized as follows.   Pick   $h\in R$ and observe that
 the  mapping $ h \nu :A\to A\otimes A$ is  a Lie cobracket in
$A$.  It induces the Lie bracket $[,]_h=h [,]$ in $A^*$ where $[,]$ is the Lie bracket induced by $\nu$.  The corresponding
multiplication
$\mu_h$  in
$A^*$ is given by
$$\mu_h (a,b)= a+b + \frac{h}{2} [a,b]+ \frac{h^2}{12} ([a,[a,b]] +[b, 
[b,a]])+ 
\cdots$$
This multiplication   makes $A^*$ into a group denoted  $\Exp_h A^*$. As above, $\mu_h$  induces a comultiplication in
the symmetric algebra $S=S(A)$. This makes $S$  into a Hopf algebra  over $R$ denoted $ S_h(A)$.  For
$h=1$,  
we obtain the same objects as in the previous paragraphs.  Note for the record that for any $h\in R$, the formula  
$a\mapsto ha: A^*\to A^*$ defines a group homomorphism $\Exp_h A^*\to \Exp  A^*$.  If $h\in R$ is a non-zero-divisor,
this homomorphism is injective.

It is clear that the construction of  $\Exp_h A^*$ and $S_h (A)$ is functorial. For a  Lie  coalgebra homomorphism
$\psi$ from  
$A$ into a spiral Lie  coalgebra 
$A'$, the dual   homomorphism   $\psi^*:(A')^*\to A^*$ preserves 
  the 
Lie bracket $[,]_h$ and the group multiplication $\mu_h$.  The  algebra   homomorphism $S_h (A)\to 
S_h (A')$ induced by $\psi $ is  a homomorphism of Hopf
algebras.

 \begin{lemma}\label{l:73} The Lie coalgebra of strings $(\A_0=\A_0(R),\nu)$ defined in
Section \ref{su71}    is spiral.
		   \end{lemma}
                     \begin{proof} This follows from
the obvious fact 
that
 $\nu^{(n)}(\langle \alpha\rangle)=0$ for any string
$\alpha$ of rank $\leq 
n$. 
(Actually a stronger assertion holds: $\nu^{(n)}(\langle
\alpha\rangle)=0$ 
for any 
string $\alpha$ of rank $\leq 4n+2$.) \end{proof}

  Applying the constructions above to    the Lie coalgebra  $\A_0 $   and any  
$h\in R$, we obtain a group   
 $\Exp_h\A_0^* $ and a Hopf algebra $ S_h( \A_0 )$ over $R$. Note that  
$ S_h( \A_0 )= R[\str_0]$ as algebras.

\begin{theor}\label{th:fge} For  $R=\QQ[z]$ and $h=z\in R$, the 
           homomorphism  $\nabla: \mathcal {E} \to R[\str_0]=S_h( \A_0 )$  is an isomorphism of Hopf algebras.
                     \end{theor}

The proof of this theorem follows  the lines of \cite{tu}, Section 12 and Lemma 13.5; we omit the details.  

  \subsection{Remarks.}\label{a21a}  1. The equality $\A=\A_0 \oplus R $ implies that 
 $\Exp  \A^*=\Exp   \A_0^* \times \underline R $ where $ \A_0^*=( \A_0)^*$ and $\underline R$ is the additive group
of 
$R$.  More generally, for any $h\in R$,  we have $\Exp_h  \A^*=\Exp_h   \A_0^* \times \underline R $.

2. For any $h\in R$, the Lie
coalgebra    
$\A_{r,g} 
$ defined in Section  \ref{su735}  gives rise to  a  group $\Exp_h\A_{r,g}^* $ and a Hopf algebra $ S_h ( \A_{r,g} )$ which
are quotients of  
$\Exp_h
\A_0^* $ and   $S_h (\A_0) $, respectively. 
  The Lie coalgebra  $Z=Z(R)$  discussed in Section \ref{su75}  is  known to
be spiral, so that we have the 
associated 
  group $\Exp_h Z^*$  and  the  Hopf
algebra $S_h (Z)$. The homomorphism     $\psi_0:Z\to \A_0$  induces
  a group 
homomorphism 
$\Exp_h\A_0^*\to \Exp_hZ^*$  and   a  Hopf
algebra homomorphism   $S_h(Z)\to 
S_h (\A_0)$. 
 
 \subsection{The algebra of   open strings.}\label{a21aee} We begin with algebraic preliminaries.  Recall  
that a {\it module}  over a Lie algebra
$(L, [ \,,]:L^{\otimes 2} \to L)$ over
$R$ can be defined as  an
$R$-module
$M$ endowed with an
$R$-linear homomorphism $\rho:L\otimes M \to M$ such that \begin{equation} \label{drf} \rho ([\,,]\otimes \id_M)= \rho
(\id_L
\otimes
\rho )(\id_{L\otimes L \otimes M}- 
 \Perm_L\otimes \id_M) :L\otimes L \otimes M \to M\end{equation}
where $\Perm_L$ is  the permutation 
$x\otimes 
y\mapsto y\otimes x$ in $L^{\otimes 2}=L\otimes  L$.  (Formula \ref{drf} is equivalent to the usual identity
$[x,y]m=x(ym)-y(xm)$ for
$x,y\in L, m\in M$.)  Dually,  a  {\it comodule}  over
a Lie coalgebra 
$(A,\nu:A\to A^{\otimes 2})$ over
$R$ is  an
$R$-module
$M$ endowed with an
$R$-linear homomorphism $\rho:M\to A\otimes M  $ such that 
\begin{equation} \label{drref} (\nu \otimes \id_M) \rho =  (\id_{A\otimes A \otimes M}-
\Perm_A\otimes \id_M) (\id_A
\otimes
\rho)  \rho:M\to A\otimes A \otimes M.\end{equation} 
Such
$M$  is automatically a module over the dual Lie algebra
$A^*=\Hom_R(M,R)$: an element  $a\in A^*$ acts on $M$  by the endomorphism $\varphi_a:M\to M$ sending  $m\in M$ to 
$ - (a
\otimes
\id_{M})
\rho(m)\in R\otimes M=M$. 

 A comodule
$ (M,\rho)$ over a Lie coalgebra $A$  is {\it spriral} if $M=\cup_{n\geq 1} \Ker \rho^{(n)}$
where 
$$\rho^{(n)}=(\id_A^{\otimes (n-1)} \otimes \rho) \circ \cdots \circ 
  (\id_A  \otimes \rho) \circ \rho:M \to A^{\otimes n}\otimes M.$$  If  both $A$ and $M$ are 
spiral, then  
 the  action of $A^*$ on $M$ integrates into a group action of the group 
$\Exp  A^*$ on  $M$ defined by   $$a m=e^{\varphi_a} (m)=m+\sum_{k\geq 1}   (\varphi_a)^k(m)/k!$$
for  $a\in \Exp  A^*=A^*$, $m\in M$. Note that for   $m\in \Ker \rho^{(n)}$ the sum on the right-hand side has at
most $n$ non-zero terms.

  Let
$\M=\M (R)$ be the free 
$R$-module freely generated by the 
set of 
homotopy classes of open virtual strings.  We    provide $\M$ with the 
structure of a comodule over the  
  Lie coalgebra $\A_0$.  Let $\langle \beta \rangle$ be the generator of $\M$ represented by an open  string $\beta$. For
an arrow $e\in \arr(\beta)$,
  a surgery along $e$ defined as in Section \ref{vks:624} transforms $\beta$ into a disjoint union of a  closed  string
$\alpha_e$ and an open string $\beta_e$.  Set  
$$\rho (\langle \beta \rangle) =\sum_{e\in \arr_+(\beta)}   \langle \alpha_e \rangle \otimes \langle \beta_e
\rangle  - \sum_{e\in \arr_-(\beta)}   \langle \alpha_e \rangle \otimes \langle \beta_e
\rangle  \in \A_0 \otimes \M.$$
A direct computation shows that this  gives a well-defined  $R$-linear homomorphism $\rho:\M\to \A_0\otimes \M 
$ satisfying Formula \ref{drref}.  Thus $\M$ is a   comodule over  $\A_0$.   
Combining $\rho$ with the inclusion $\A_0 \subset \A$ we obtain that $\M$  is a
comodule over   $\A$ as well.   It is easy to see that $\M$ is spiral.  The construction
 above gives a group action of
$  \Exp  \A^*$ on  $\M$.  

\subsection{Exercises.}\label{a21ass}  1. The obvious multiplication of open strings makes $\M$ into an associative algebra
with unit. Check  that  the group $\Exp   \A^*$ acts on  $\M$ by algebra automorphisms. 

2.   Let $cl:\M\to \A$  be the   $R$-linear homomorphism induced by closing open strings.  Check that for any open
string
$\beta$,
 we have $\nu(\langle cl(\beta)  \rangle)=(\id_{\A\otimes \A}  -\Perm_{\A} )  (\id_\A \otimes cl) \rho
(\langle \beta \rangle)  $.

	\section{Open questions}

    1.  Find out whether the slice genus of  a string is a homotopy invariant. 
The slice genus is invariant under  
(c$)_s$ but may possibly decrease under  (a$)_s$, 
(b$)_s$. Find out whether   these moves can transform a non-slice string into a slice one.

2.  Find 
out which  primitive based matrices  $T_0$ can  be realized as $T_0(\alpha)$ for a string $\alpha$.     A
necessary  condition 
pointed out in Section \ref{sn:g22}  says that $u'_{T_0}(1)=0$. 
Are there other conditions ~?  Note that for the based matrix 
$T(\alpha)=(G,s,b)$, we have $\vert b(e,f)\vert \leq \# 
(G)-2$ for 
all $e,f \in G$.  
     This   however yields no conditions on the primitive based matrices arising from strings, since such a matrix 
$T_0=(G_0,s,b_0)$ may  arise from a 
string of a rank much bigger than $\# (G_0)$.

3. Find further obstructions  to the sliceness of a string.
Specifically, are the virtual strings  $\alpha_{p,p}$ with $p\geq 2$ slice (cf. Corollary \ref{th:t254711}) ? The based
 matrix 
 of $\alpha_{p,p}$ is hyperbolic  and   gives no information on the question.

4.   Consider  the string of rank four
$\alpha=\alpha_\sigma$ where
  $\sigma=(1342)$. A direct computation shows  that its primitive based matrix
 is trivial. Also $\nu (\langle \alpha \rangle)=0$  since $\alpha $ has only 4 arrows. Is $\alpha$
homotopically trivial  ? A more ambitious program would be to classify all strings of small rank (say, $\leq 6$) up to
homotopy.

5. Is  multiplication of open strings discussed in Section \ref{sn:g14}  commutative (up to homotopy) ? In other words, is
the algebra of   open strings  $\M$  considered at the end of  Section \ref{fur}   commutative?

6.    Is there an  invariant of formal knots combining the skein invariant
$\nabla$   with the Kontsevich universal finite type invariant of knots ?  This might lead to  mixed arrow-chord
diagrams. 

7.  Study   invariants of virtual strings that change 
in a 
controlled way (say by constants) under the moves (a$)_s$, (b$)_s$,  (c$)_s$, 
cf. the 
theory of Arnold's invariants of plane curves.

                     \end{document}